# TRIMMING AND LIKELIHOOD: ROBUST LOCATION AND DISPERSION ESTIMATION IN THE ELLIPTICAL MODEL[1]

By Juan A. Cuesta-Albertos, Carlos Matrán and Agustín Mayo-Iscar

*Universidad de Cantabria, Universidad de Valladolid and Universidad de Valladolid*

Robust estimators of location and dispersion are often used in the elliptical model to obtain an uncontaminated and highly representative subsample by trimming the data outside an ellipsoid based in the associated Mahalanobis distance. Here we analyze some one (or $k$)-step Maximum Likelihood Estimators computed on a subsample obtained with such a procedure.

We introduce different models which arise naturally from the ways in which the discarded data can be treated, leading to truncated or censored likelihoods, as well as to a likelihood based on an only outliers gross errors model. Results on existence, uniqueness, robustness and asymptotic properties of the proposed estimators are included. A remarkable fact is that the proposed estimators generally keep the breakdown point of the initial (robust) estimators, but they could improve the rate of convergence of the initial estimator because our estimators always converge at rate $n^{1/2}$, independently of the rate of convergence of the initial estimator.

## 1. Introduction.

Between the methodologies to produce robust and efficient estimators we are here concerned with those based on a preliminary robust estimation followed by one step (or $k$ steps) that improves efficiency without a significant loss of robustness. In a natural way this leads us to search for an uncontaminated and highly representative subsample, selected

Received July 2007; revised July 2007.

[1]Supported in part by the Spanish Ministerio de Educación y Ciencia and FEDER, Grant BFM2005-04430-C02-01, and 02 and by the Consejería de Educación y Cultura de la Junta de Castilla y León, Grant PAPIJCL VA102A06.

using the initial estimation, and then to make the improvement step on the basis of this subsample. These ideas are present, for example, in Rousseeuw and van Zommeren [25] or in Lopuhaä and Rousseeuw [19], where it is shown that some ways of one-step reweighting preserve the breakdown point (BP) of the initial estimators. This scheme seems to be particularly adequate to preventing gross errors under a model based on a main data stream of an elliptical distribution. The robust estimates are then used as a diagnostic tool to select the good observations, such as those at an adequate (Mahalanobis) distance from the location estimate. Hence, we could improve the efficiency of our estimators, preserving robustness with respect to outliers, by resorting to classical methods which obtain efficient estimates but compute only over the observations considered good.

In this brief description three main ingredients require consideration:

- The choice of the robust initial estimator to produce the zone of good observations, that is, a suitably trimmed set.
- Once such a zone has been selected, how to treat the discarded (trimmed) data.
- How to choose the efficient estimator.

The first item has received considerable attention in the elliptical model, which allows us to exploit the symmetries in order to handle gross-errors as points far away from the center. Under the usual equivariance requirement these robust estimators include well-known proposals like the Minimum Volume Ellipsoid (MVE), the Minimum Covariance Determinant (MCD) or, in general, S-estimators (see the book by Maronna, Martin and Yohai [22] for a discussion of these and further estimators in this setup).

The other items have been treated in an unequal way. Usually the one-step consists of reweighted least squares statistics based only on good sample data, which take advantage of the elliptical symmetry of the underlying family (see, e.g., Lopuhaä [18]). Some versions, as in Gervini [12], exploit even the possibility of selecting the good sample data region in an adaptive way. A different point of view (see, e.g., Bickel [2] or Davies [8]) resorts to a Newton–Raphson step to increase the rate of convergence of the initial estimators. Curiously, the maximum likelihood estimator (MLE), being a natural choice in order to get the maximum gain in efficiency, has only recently been considered in Mayo-Iscar [23] in the mixture model, and in Marazzi and Yohai [20] in the context of regression of a real valued variable. However, under simple truncation (which is the approach followed in [20]), existence of the corresponding MLE in any data configuration is not guaranteed, so breakdown of the estimator could arise under contamination leading to such a configuration.

Here we will mainly address the last two stressed items. The starting point will be that of a given (trimmed) set of the sample space, obtained through any of the already enumerated methods and that is consistent under the subjacent model. The consistency of the mentioned methods is well known (see, e.g., Davies [6, 7] and Butler, Davies and Jhun [4]) under the elliptical model, but has not been treated in a contamination model like our GEM proposal below. Since the MVE is better adapted to this contaminated model and it is the more impacting possibility, even in the uncontaminated model because of its low efficiency, this estimate will act as a leitmotif in the paper. Its consistency is shown in Proposition 3.1.

With respect to the second item, first we will consider the likelihoods associated to the (artificially) truncated and censored models given the set, but we will also introduce a model of gross errors contamination and consider the associated likelihood. In every case, in connection with the third item, we will consider the MLEs. The relations between the MLE associated to each likelihood model provide a novel approach to interpret the presence of gross errors. Under our model we could adapt the final estimation through a second step, based on a cut-off parameter, in a similar way to that introduced in Gervini and Yohai [13] (improved in [20]) or in García–Escudero and Gordaliza [11]. However, in order to make the comprehension of the methodology easier, we will not consider here these adaptive ways of enlarging the MVE to improve the final efficiency of the estimator, although the corresponding analysis would be parallel to the one developed here.

An important feature of our approach is the rate of convergence of the final estimator. As a distinctive fact with respect to the known one-step reweighted estimators, our estimators converge at rate $n^{1/2}$, independently of the rate of convergence of the initial estimator, whenever it is consistent. This allows any consistent initial estimator to be considered, even the MVE that converges at $n^{1/3}$ rate, as an initial estimator without loss in the rate of convergence. This happens because, although based on only a part of the sample, our second step is a genuine MLE, so it is able to make a full reconsideration of the initial estimation. On the contrary, this is not possible if we only make a linear estimation based on reweighting in accordance with the initial estimation. These considerations agree with those in Rousseeuw [24] or He and Portnoy [14], where it is stressed that the problem is reweighting. However, the estimators considered in the already mentioned literature as one-step improvements to get the $n^{1/2}$-rate of convergence are based on a Newton–Raphson adjustment (see also Jurečkova and Portnoy [16] or Jurečkova and Sen [17]). In fact, in [20], it is even suggested (see Remark 2 after Theorem 3 there) that the rate of convergence of the truncated MLE could be the same as that of the initial estimator.

The paper is organized as follows. In Section 2 we introduce the models and the estimators to be studied, and analyze the identifiability of the models. Moreover, we discuss the existence and uniqueness of the MLEs under the truncated, the censored and the gross errors models (GEM). We stress the fact that truncation can produce the nonexistence of the MLE, and hence produce the breakdown of the estimator; thus, we introduce natural restrictions which guarantee the existence of the MLE under truncation. In this setting we obtain new results, even in the univariate case, which include the consideration of exponential families. Section 3 is devoted to studying the robustness and the asymptotic properties of the proposed estimators, including the BP and the influence function (IF), as well as the consistency and asymptotic normality at rate $n^{1/2}$. In Section 4 we present our conclusions on these estimators. The paper ends with an Appendix containing all the proofs and some technical results.

The estimation will be carried out on the basis of the data set $\mathcal{X} = \{x_1, x_2, \ldots, x_n\}$. In the asymptotic results we will assume that $\mathcal{X}$ is obtained from $n$ independent, identically distributed $\mathbb{R}^p$-valued random vectors $X_1, X_2, \ldots, X_n$. $P_n$ will be the associated sample distribution. The usual norm in $\mathbb{R}^p$ will be denoted by $\| - \|$, the $\sigma$-field of Borel sets will be $\beta^p$ and $\lambda^p$ will be the Lebesgue measure. When we use matrix notation, vectors must be understood as column vectors. For a matrix $H = (h_{ij})$, $H^T$ will denote its transpose and $|H|$ its determinant. $B(m,r)$ [resp. $S(m,r)$] will denote the open ball (resp. sphere) of radius $r$ centered at $m$, while $A^c$ is the complement of the set $A$, and $I_A$ is its indicator function. Further notation will be introduced throughout the paper as necessary.

We will make use of a generic $\omega$ in a probabilistic space of reference; almost surely (a.s.) statements must be understood as relative to that space. $P_n^\omega$ would be then a realization of $P_n$. Integration of a random variable $h$ with respect to a probability $\mathbb{P}$ will be denoted as $\mathbb{P}h$ (and is interpreted componentwise when $h$ is a vector). We will use matrix notation for partial derivatives of a function. Given $g : \Theta \times \Phi \times \mathbb{R}^p \to \mathbb{R}$, where $\Theta \subset \mathbb{R}^d, \Phi \subset \mathbb{R}^k$, $\frac{\partial}{\partial\theta} g(\theta, \phi, x)$ will denote the $d$-dimensional vector with components $\frac{\partial}{\partial\theta^i} g(\theta, \phi, x), i = 1, \ldots, d$.

In our setup the choice of an initial set leads us to consider ellipsoids indexed by the set $\Gamma := \mathbb{R}^p \times \mathcal{M}^+_{p\times p} \times \mathbb{R}^+$, where $\mathcal{M}^+_{p\times p}$ is the set of positive-definite symmetric $p \times p$ matrices. For $\gamma = (\mu, \Sigma, r) \in \Gamma$, we will denote $\mathcal{E}(\gamma) := \{x \in \mathbb{R}^p : (x-\mu)^T \Sigma^{-1} (x-\mu) \le r^2\}$.

**2. Maximum likelihood estimation with a trimmed sample.** We begin with the minimal assumptions that we consider throughout this paper.

Definition 2.1. The elliptical model associated to the nonincreasing function $g : \mathbb{R}^+ \to \mathbb{R}^+$ is a family $\{\mathbb{P}_\theta : \theta \in \Theta\}$ of probabilities on $\beta^p$ with densities $f_\theta$ (with respect to $\lambda^p$) given by

$$f_\theta(x) = |\Sigma|^{-1/2} g((x-\mu)^T \Sigma^{-1} (x-\mu)), \tag{2.1}$$

where $\Theta := \mathbb{R}^p \times \mathcal{M}^+_{p\times p}$, and $\theta := (\mu, \Sigma) \in \Theta$. Note that $\mathcal{M}^+_{p\times p}$, considered as a subset of $\mathbb{R}^{\binom{p+1}{2}}$, is open, $\mu$ is the mean of $\mathbb{P}_\theta$ when it exists, while, if the second moment is finite, the variance-covariance matrix of $\mathbb{P}_\theta$ is proportional to $\Sigma$.

We will also handle a contaminated version of this model: The *only gross-errors outliers*. This term is used to indicate that $\mathbb{P}_\theta$ could be contaminated with a small proportion of data coming from a distribution whose support is external to a central part of $\mathbb{P}_\theta$. Throughout, a *central part of an elliptical distribution* must be understood as an ellipsoid $\mathcal{E}(\mu, \Sigma, r)$ where $\mu$ and $\Sigma$ are the location and scatter parameters of the distribution. These sets have the nice property of being scaled versions of the MVE; see Lemma A.3 (but also of the MCD, see Butler, Davies and Jhun [4]). Recall that *an ellipsoid $A$ is an MVE of $\mathbb{P}$, if $\mathbb{P}(A) \geq 1/2$ and the volume of $A$ is minimal in the class of all ellipsoids with this property.*

Existence and uniqueness of MVE's are discussed in Davies [7]. As stated, if $\mathbb{P}_{(\mu,\Sigma)}$ belongs to the elliptical model, then their MVE's are essentially given by ellipsoids $\mathcal{E}(\mu, \Sigma, r)$, where $r$ depends on $\mu$ and $\Sigma$ but also on $g$ and $p$. In particular, when $g$ is strictly decreasing the MVE is unique. Below we describe our version of the GEM that is considered from now on.

DEFINITION 2.2 (Gross error model). A distribution $\mathbb{P}$ belongs to the Gross Error Model associated to the family $\{\mathbb{P}_\theta : \theta \in \Theta\}$ if there exist $\pi \in [0, 1)$, a probability $Q$ and $\theta \in \Theta$, such that

$$\mathbb{P} = (1 - \pi)\mathbb{P}_\theta + \pi Q, \tag{2.2}$$

where $Q$ is a probability distribution such that if $A$ is any MVE of $\mathbb{P}$, then $Q(A) = 0$ (whence $A$ is also a central part of $\mathbb{P}_\theta$).

Our proposal to produce the estimator is this: First, through a consistent estimator of $\theta$, we produce an estimation $\mathcal{E}(\hat{\mu}, \hat{\Sigma}, \hat{r})$ of the MVE of $\mathbb{P}$ because, even in the GEM, at least asymptotically, the values in the sample that remain in the estimated MVE would be produced by the elliptical part of $\mathbb{P}$. Then, in order to maintain the BP of the initial estimator (see Theorem 3.1) while achieving the highest possible efficiency, we will enlarge the estimated MVE by keeping their location, $\hat{\mu}$, and shape, $\hat{\Sigma}$, but taking a greater value than $\hat{r}$ to get a scaled (thus containing more points of the sample) version of the estimated MVE. As a final step we will construct an MLE of $\mu$ and $\Sigma$ based on the observations lying in this scaled MVE.

This program essentially coincides with that introduced in [20], although, to make the exposition easier, we will focus on the one-step estimator based on the (nonenlarged) MVE. However, in Tables 1, 2 and 3, in Appendix B, we will present the gains in efficiency attained by handling the scaled versions.

2.1. *Identifiability of the model.* Since we will discard the data not included in an ellipsoid, we need to assure that the parameters of an elliptical distribution are identifiable from every ellipsoid or, more generally, from every open set. This holds in the Gaussian case and other general models like the exponential family, as we show in Proposition 2.1.

PROPOSITION 2.1. *Let $\{f_\theta : \theta \in \Theta\}$ be the density functions of a $d$-parameter exponential family with respect to a $\sigma$-finite measure $\lambda$ on $\mathbb{R}^p$, where*

$$f_\theta(x) = C(\theta)\exp\left\{\sum_{j=1}^{d} Q_j(\theta)T_j(x)\right\}h(x). \tag{2.3}$$

*Assume that the $Q_j$'s do not satisfy a linear constraint on $\Phi$. Let $A \in \beta^p$ such that the $T$'s $\lambda$-a.s. do not satisfy a linear constraint on $A$.*

*If $f_{\theta_1} = f_{\theta_2}$ $\lambda$-a.s. on $A$, then $\theta_1 = \theta_2$.*

The identifiability shown in Proposition 2.1 can be circumvented because we only need to guarantee that each probability $\mathbb{P}_\theta$ is identifiable from the rest over adequate sets. Broadly speaking, we can say that a set $A$ is adequate for $\mathbb{P}_\theta$ if $\mu \in A$ and $f_\theta$ is not constant on $A$.

Since $g$ is nonincreasing, to assure that $g$ is neither constant on $A$ nor on affine transformations of $A$, a natural hypothesis is to assume that $g$ is strictly decreasing. But, since the ellipsoids of interest (asymptotically) contain $\mu$, it is enough to demand that $g$ is strictly decreasing near zero, or, more generally that $g$ fulfills the following condition:

(G1) There exists a strictly decreasing sequence $\{t_n\}$ which converges to zero, such that $g(t_n) < g(t_{n+1})$ for every $n$.

PROPOSITION 2.2. *Let $\{\mathbb{P}_\theta : \theta \in \Theta\}$ be the elliptical model on $\mathbb{R}^p$, $p \geq 1$ associated to $g$ which verifies condition G1. Let $\theta_0 = (\mu_0, \Sigma_0) \in \Theta$ and $A$ be an open set in $\mathbb{R}^p$, such that $\mu_0 \in A$. If $\theta \in \Theta$, $\theta \neq \theta_0$, then*

$$\mathbb{P}_{\theta_0}\{x \in A : f_\theta(x) \neq k f_{\theta_0}(x)\} > 0, \qquad \textit{for every } k > 0.$$

2.2. *The estimators.* Once we know that it is possible to estimate the elliptical part from adequate sets, we will analyze some estimation procedures related to this task. In Section 2.2.1 we will consider the MLE associated to (in our case artificially) truncated or censored samples. In Section 2.2.2, under the GEM, we design a new estimator, called the Smart estimator. In every case their effective computation can be implemented through the EM algorithm. Section 2.2.3 explores the existence and uniqueness of these estimators.

2.2.1. *Censored and truncated estimators.* The difference between truncated and censored estimators lies on the way in which they consider the discarded points in the sample. The censored one forgets the right values of the points outside $A$, but takes into account their number. Thus, this number should appear in the likelihood function but related to no specific point. To this, we introduce an artificial point $\mathbf{c}$, not necessarily in $\mathbb{R}^p$, which is only used to count the number of censored points. Thus, the objective function is the *censored log-likelihood*:

$$L^c_{\theta/A}(x) := I_A(x)\log f_\theta(x) + I_{\{\mathbf{c}\}}(x)\log \mathbb{P}_\theta(A^c), \qquad x \in \mathbb{R}^p \cup \{\mathbf{c}\}; \tag{2.4}$$

here the points in $A^c$ are treated as located on the identical censored state $\mathbf{c} \notin A$, and $0 \times \infty$ is taken as 0. This log-likelihood corresponds to the model $\{\mathbb{P}^c_{\theta,A} : \theta \in \Theta\}$ given by $\mathbb{P}^c_{\theta,A}(B) = \mathbb{P}_\theta(B \cap A) + \mathbb{P}_\theta(A^c)I_B(\mathbf{c}), B \in \sigma(\beta^p \cup \{\mathbf{c}\})$, with density function $f^c_{\theta,A}(x) = f_\theta(x)I_A(x) + \mathbb{P}_\theta(A^c)I_{\{\mathbf{c}\}}(x), x \in \mathbb{R}^p \cup \{\mathbf{c}\}$, with respect to the measure $\lambda^p + \delta_{\{\mathbf{c}\}}$, where $\delta_{\{\mathbf{c}\}}$ is the Dirac measure on $\mathbf{c}$.

Obviously the MLE based on $L^c_{\theta/A}$ is more stable than the usual MLE in presence of some contamination in $A^c$. However, this contamination can produce an excessive weight on $A^c$. To protect against this possibility, we can consider truncation, which does not consider at all the data in $A^c$, the *truncated log-likelihood* being

$$L^t_{\theta/A}(x) := I_A(x)\log\frac{f_\theta(x)}{\mathbb{P}_\theta(A)}, \qquad x \in \mathbb{R}^p. \tag{2.5}$$

This corresponds to the model $\{\mathbb{P}^t_{\theta,A} : \theta \in \Theta\}$ defined through the density functions $f^t_{\theta,A}(x) = I_A(x)f_\theta(x)/\mathbb{P}_\theta(A), x \in \mathbb{R}^p$ with respect to $\lambda^p$. In agreement with the obvious incompatibility which would arise for those $\theta$'s such that $\mathbb{P}_\theta(A) = 0$, we adopt the convention that $L^t_{\theta/A} = -\infty$ if this happens.

Thus, given a sample $\mathcal{X} = \{x_1, \dots, x_n\}$, the *maximum likelihood censored* (*resp. truncated*) *estimator*, MLE(c) [resp. MLE(t)] on an appropriately chosen set $A$ will be the value $\hat{\theta}_{c,n}$ (resp. $\hat{\theta}_{t,n}$) maximizing $P_n L^c_{\theta/A}$ (resp. $P_n L^t_{\theta/A}$).

2.2.2. *A data based choice: The smart estimator.* We present an estimator which takes full advantage of the GEM. As far as we know, likelihood-based estimators under this model have not yet been proposed.

We have to face two difficulties: There is not a unique set $A$ related to the model and, given an observation in $A^c$, we do not know whether this observation comes from $\mathbb{P}_\theta$ or from $Q$. To circumvent the first difficulty, given a sample of size $n$, for every suitable set $A$, we can consider the log-likelihood of $n_1$ (resp. $n_2$) data points in $A$ (resp. in $A^c$) arising from $\mathbb{P}_\theta$ and $n_3 = n - n_1 - n_2$ from the contaminating source also in $A^c$. For the second, we consider a model in which we have the complete information for the data in $A$ and only the global $(n_2+n_3)$ number of points in $A^c$.

Thus, we have to make a first estimation to get a suitable set $A$, which can be understood as a noise parameter in the model, and realize the final estimation on the basis of the likelihood associated to this empirical set. In analogy with the censored likelihood in (2.4), we consider an ideal censored state $\mathbf{c}$ and, for every $x \in \mathbb{R}^p \cup \{\mathbf{c}\}$, define the log-likelihood

$$\begin{aligned} L^s_{\theta,\pi/A}(x) := {} & I_A(x)\log((1-\pi)f_\theta(x)) \\ & + I_{\{\mathbf{c}\}}(x)\log((1-\pi)\mathbb{P}_\theta(A^c)+\pi), \end{aligned} \tag{2.6}$$

which is associated to the model $\{\mathbb{P}^s_{\theta,\pi,A} : \theta \in \Theta, \pi \in [0,1)\}$ given by

$$\mathbb{P}^s_{\theta,\pi,A}(B) = (1-\pi)\mathbb{P}_\theta(B\cap A) + ((1-\pi)\mathbb{P}_\theta(A^c)+\pi)I_B(\mathbf{c}), \tag{2.7}$$

where $B \in \sigma(\beta^p \cup \{\mathbf{c}\})$. The sample objective function to maximize is now

$$P_n L^s_{\theta,\pi/A} = P_n(I_A\log((1-\pi)f_\theta) + I_{A^c}\log((1-\pi)\mathbb{P}_\theta(A^c)+\pi)), \tag{2.8}$$

under the restriction $\pi \in [0,1)$. The estimator obtained by maximizing this objective function will be called *smart MLE* [MLE(s)]. The analysis of the existence of this estimator, to be carried in the next subsection, will shed new light on our proposal for this problem.

2.2.3. *On the existence and uniqueness of the estimators.* The existence and uniqueness of the MLE is not an easy problem. In fact, the truncated normal model is often used as an example of possible inexistence of the MLE.

For the elliptical model, Maronna [21] treated the problem of existence and uniqueness of M-estimators, for the model and the sample, but his assumptions on $g$ are not satisfied, for example, by the normal model or by our models related to truncation or censoring.

Under the theoretical model both facts are an easy consequence of Jensen's (strict) inequality and the identifiability. The proof is similar to the classic one.

PROPOSITION 2.3. *Under the hypotheses of Proposition 2.2, for every $\theta \neq \theta_0$, we have*

$$\mathbb{P}^c_{\theta_0,A} L^c_{\theta_0/A} > \mathbb{P}^c_{\theta_0,A} L^c_{\theta/A} \tag{2.9}$$

*and*

$$\mathbb{P}^t_{\theta_0,A} L^t_{\theta_0/A} > \mathbb{P}^t_{\theta_0,A} L^t_{\theta/A}, \tag{2.10}$$

*and, for every $(\theta,\pi) \in \Theta \times [0,1) - \{(\theta_0,\pi_0)\}$,*

$$\mathbb{P}^s_{\theta_0,\pi_0,A} L^s_{\theta_0,\pi_0/A} > \mathbb{P}^s_{\theta_0,\pi_0,A} L^s_{\theta,\pi/A}. \tag{2.11}$$

To obtain the existence of the MLEs, we should avoid, for a sample in general position, a degenerated (into a lower dimension) solution. This is related to the speed of decreasing of $g$ and leads us to introduce the following assumption Gp. Moreover, we will impose the continuity of $g$ as another natural requirement:

(Gp) If $p > 1$, then there exists $\gamma > p/2$ such that $\lim_{r\to\infty} r^\gamma g(r) = 0$.
(G2) $g$ is continuous on $\mathbb{R}^+$.

Note that, by Scheffé's lemma, G2 implies that

$$\sup_{A\in\beta^p} |\mathbb{P}_{\theta_n}(A) - \mathbb{P}_{\theta_0}(A)| \to 0, \qquad \text{whenever } \theta_n \to \theta_0. \tag{2.12}$$

Proposition 2.4 (Existence of nonrestricted MLE). *Let $g$ be a function which defines an elliptical family on $\mathbb{R}^p$ and satisfies G1, G2 and Gp. Let $n > 2$, if $p = 1$, and let $n > \frac{p\gamma}{\gamma - p/2}$ in the case $p > 1$, where $\gamma$ is the constant which appears in Gp.*

*Then, for every data set $\mathcal{X} = \{x_1, \dots, x_n\}$, whose points are in general position, there exists $(\hat\mu, \hat\Sigma) \in \Theta$ such that*

$$\prod_{i=1}^n f_{(\hat\mu,\hat\Sigma)}(x_i) \ge \prod_{i=1}^n f_{(\mu,\Sigma)}(x_i), \qquad \textit{for every } (\mu,\Sigma) \in \Theta.$$

The next proposition proves the existence of the smart and censored estimators.

Proposition 2.5 [Existence of MLE(s) and MLE(c)]. *Assume that $g$ defines an elliptical family on $\mathbb{R}^p$ and satisfies G1, G2 and Gp, and let $\mathcal{X} = \{x_1, \dots, x_n\}$ be a data set.*

*Let $A \in \beta^p$ such that the number of points, $m$, in the set $\mathcal{X} \cap A$ satisfies that $m > 2$, if $p = 1$, and that $m > \frac{p\gamma}{\gamma - p/2}$ in the case $p > 1$, where $\gamma$ is the constant which appears in Gp.*

*If the points in $\mathcal{X} \cap A$ are in general position, then there exist the MLE(s), $(\hat\theta_{s,n}, \hat\pi_n)$, and the MLE(c), $\hat\theta_{c,n}$, based on the sample $\mathcal{X}$ and $A$.*

The existence of the MLE(t) cannot be shown with the same argument because the denominator $\mathbb{P}_{\theta_k}(A)$ in (2.5) could converge to zero. In fact, this can lead to nonexistence of the MLE(t). This difficulty can be handled on the basis that the sets $A$ under consideration will be estimations of the MVE of $\mathbb{P}$, thus their probabilities must be large enough.

Given $\alpha > 0$ and the ellipsoid $A$, let

$$\Theta_A^\alpha := \{\theta \in \Theta : \mathbb{P}_\theta(A) \ge \alpha\}. \tag{2.13}$$

Assume that $\mathbb{P}$ is a probability obtained by contaminating $\mathbb{P}_{\theta_0}$ by any probability $Q$ with $Q(A^c)=1$, and $\mathbb{P}_{\theta_0}(A)=\alpha_0$. We will obtain in Proposition 3.2 that, asymptotically, if $\alpha \in (0,\alpha_0)$, the restrictions $\Theta^{\alpha}_{A_n}$ obtained from the sample MVE, $A_n$, are satisfied by every $\theta$ in a neighborhood of $\theta_0$. Moreover, as stated below, the truncated likelihood function constrained to the set $\Theta^{\alpha}_{A_n}$ has a maximum. These facts allow us to consider the MLE(r) or *constrained MLE*(*t*), to be denoted as $\hat{\theta}_{r,n}$, as a substitute of the MLE(t).

PROPOSITION 2.6. *Given $\alpha>0$, let $\Theta^{\alpha}_A$ be defined as in (2.13). Let us assume the hypotheses in Proposition 2.5 for $g$, $\mathcal{X}$ and $A$.*

*If the points in $\mathcal{X}\cap A$ are in general position, then there exists $\hat{\theta}_{r,n}\in\Theta^{\alpha}_A$, such that*

$$P_n L^t_{\hat{\theta}_{r,n}/A} = \sup_{\theta\in\Theta^{\alpha}_A} P_n L^t_{\theta/A}.$$

Dependence of the constrained solutions on the $\alpha$-value could be considered as a drawback of this proposal. However, as shown in the next proposition, in our setup the level defining the restriction will arise in a natural way, justifying our considering the MLE(r) for $\alpha = P_n(A)$ as a *natural* MLE(r).

PROPOSITION 2.7. *Let us assume the hypotheses in Proposition 2.5.*

*Let $(\hat{\theta}_{s,n},\hat{\pi}_{s,n})$ be an MLE(s) and let us define, for every $\theta\in\Theta$,*

$$\pi^*(\theta) := \frac{\mathbb{P}_\theta(A) - P_n(A)}{\mathbb{P}_\theta(A)}. \tag{2.14}$$

*If $\hat{\pi}_{s,n}=0$, then $\hat{\theta}_{s,n}$ is an MLE(c).*

*If $\hat{\pi}_{s,n}>0$, then $\hat{\pi}_{s,n}=\pi^*(\hat{\theta}_{s,n})$ and $\hat{\theta}_{s,n}$ is an MLE(r) restricted to $\Theta^{\alpha}_A$, for $\alpha = P_n(A)$.*

The key to compare the proposed estimators is the MLE(t) when it exists (see Theorem 2.1 and Proposition 3.2). From the arguments in the proof of Proposition 2.7, if $\pi^*(\hat{\theta}_{t,n})\geq 0$ then $(\hat{\theta}_{t,n},\pi^*(\hat{\theta}_{t,n}))$ would be the MLE(s), while if $\pi^*(\hat{\theta}_{t,n})<0$, then the maximum of $P_n L^s_{\hat{\theta}_{t,n},\pi/A}$ on $[0,1]$ is obtained for $\pi=0$, so the solution given by the MLE(c) and $\pi=0$ would be preferable. In other words, in spite of the MLE(c) always existing, under the assumptions in Proposition 2.5, it is only preferred when the MLE(t) produces troubles, either because the MLE(t) does not exist or because the associated estimation of $\pi$ (given by $\hat{\pi}_{t,n}$) is negative. But the MLE(t) only takes into account the data inside $A$, thus the troubles appear either because they are not likely enough to arise from the elliptical distribution, or because this estimation leads us waiting on more sample data outside $A$. Proposition 2.8 highlights these facts.

PROPOSITION 2.8. *Assume the assumptions and notation of Proposition 2.7 and that there exists an MLE(t), $\hat{\theta}_{t,n}$. Then $\pi^*(\hat{\theta}_{t,n}) \geq 0$ implies that $(\hat{\theta}_{t,n}, \pi^*(\hat{\theta}_{t,n}))$ is an MLE(s). Otherwise $(\hat{\theta}_{c,n}, 0)$ would be an MLE(s).*

It was precisely this behavior that led us to give the name "smart" to our estimate, in order to stress this suggestive property of choosing between two estimators. Whenever we make reference to the global problem, including the estimation of the contamination level, we will also use *smart estimate* to refer to the pair $(\hat{\theta}_{s,n}, \hat{\pi}_n)$, where $\hat{\pi}_n$ is defined as $\pi^*(\hat{\theta}_{s,n})$ in (2.14) when it is feasible and as 0 otherwise.

We also stress that under the GEM the consistency of the MLE(s) will imply that $\hat{\pi}_{t,n}$ is positive for large $n$, so that $(\hat{\theta}_{t,n}, \pi^*(\hat{\theta}_{t,n}))$ will asymptotically produce the smart estimator.

Uniqueness of the MLEs in our different schemes is a very distinct task. In any case, it should be noticed that the uniqueness of the estimators themselves is not necessary to obtain results on their asymptotic behavior or even their BP. In general, the treatment of the uniqueness of the MLE is closely related to the exponential family (see [1]), and this is also our approach in Theorem 2.1.

THEOREM 2.1. *Let $\{f_\theta : \theta \in \Theta\}$ be the density functions of a $d$-parameter exponential family with respect to a $\sigma$-finite measure $\lambda$ on $\mathbb{R}^p$ given by (2.3).*

*Let $A \in \beta^p$ and $\mathbb{P}$ be any probability on $\mathbb{R}^p$ such that $\mathbb{P}(A) > 0$ and $\mathbb{P}(|T_j I_A|) < \infty$, $j = 1, \dots, d$. Assume that neither the $T$'s on $A$ ($\mathbb{P}$-a.s.), nor the $Q$'s on $\Theta$ satisfy a linear constraint and let $\mathbb{P}L^s_{\theta,\pi/A}$ be the expected log-likelihood, under $\mathbb{P}$ of (2.6).*

*Then, there exists at most one solution for the maximization of $\mathbb{P}L^s_{\theta,\pi/A}$ under the restriction $\pi \geq 0$ and there exists at most one solution for the maximization of $\mathbb{P}L^t_{\theta/A}$. Moreover, if there exists a solution constrained to $\Theta^\alpha_A$, $\hat{\theta}$, which verifies $\mathbb{P}_{\hat{\theta}}(A) > \alpha$ (i.e., it is not in the boundary of $\Theta^\alpha_A$), then it is unique and also solves the unconstrained problem.*

As a consequence of Theorem 2.1 and Proposition 3.2, we can assure that the MLE(t) exists asymptotically and that it is unique for the exponential family. The following corollary particularizes this for the normal family.

COROLLARY 2.1. *Let $\{\mathbb{P}_\theta : \theta \in \Theta\}$ be the normal $p$-dimensional family, and let $A$ be any bounded set whose interior is nonempty. If $\mathcal{X} = \{x_1, \dots, x_n\}$ is a data set such that $\mathcal{X} \cap A$ has at least $p+1$ points which are in general position, then there exists a unique smart estimator $(\hat{\pi}_n, \hat{\theta}_{s,n})$, at least there exists one MLE(c), and at most there exists one MLE(t) based on $A$. Moreover, for every $\alpha \in (0,1)$ there exists an MLE(r). In particular, there exists a natural MLE(r) [corresponding to $\alpha = \mathbb{P}_n(A)$].*

2.3. *Information matrices.* This section ends with the computation of the information matrices of the proposed estimators. Those results, under the hypothesis of regularity of the model, will be employed in Theorem 3.3 to obtain the asymptotic distributions.

Regularity of a statistical experiment demands the following (see, e.g., page 65 in [15]): (a) continuity of the densities $f_\theta(x)$ on $\Theta$ for $\lambda^p$-a.e. $x$; (b) Fisher's finite information at every $\theta \in \Theta$ [i.e., differentiability of the function $f_\theta^{1/2}(\cdot)$ in $L_2(\lambda^p)$ at every point $\theta \in \Theta$], and (c) continuity in the space $L_2(\lambda^p)$ of this differential function for every $\theta \in \Theta$.

In order to guarantee the regularity of the elliptical model, we could resort to the minimal conditions given by Bickel (see pages 96–98 in [3]), consisting in the absolute continuity of $g$ and the finiteness of the integral

$$\int_0^\infty r^{p+1}(1+r^2)\left(\frac{g'}{g}\right)^2(r^2)g(r^2)\,dr.$$

Under the regularity of the statistical experiment, Lemma 7.2 in [15] shows that for any function $T$, such that $\mathbb{P}_\theta T^2$ is bounded in a neighborhood, $V_{\theta_0}$, of $\theta_0 \in \Theta$ the function $\theta :\to \mathbb{P}_\theta(T)$ is continuously differentiable in $V_{\theta_0}$ and

$$\mathbb{P}_\theta\left(T\left(\frac{\partial}{\partial\theta}\log(f_\theta)\right)\right) = \frac{\partial}{\partial\theta}\mathbb{P}_\theta(T). \tag{2.15}$$

In particular, we have $\mathbb{P}_\theta(\frac{\partial}{\partial\theta}\log(f_\theta)) = 0$.

These relations and easy computations (we omit), which take into account facts as

$$\frac{\partial}{\partial\theta}\log\mathbb{P}_\theta(A) = \frac{\frac{\partial}{\partial\theta}\mathbb{P}_\theta(A)}{\mathbb{P}_\theta(A)} = \frac{1}{\mathbb{P}_\theta(A)}\mathbb{P}_\theta\left(I_A\frac{\frac{\partial}{\partial\theta}f_\theta}{f_\theta}\right),$$

lead to the following propositions on the information matrices of our models. Notice that (except in Proposition 2.10) the involved results do not depend on the elliptical hypothesis. Proposition 2.9 also relates the information matrices based on the original, the censored and the truncated models, which we respectively denote by $\mathcal{I}(\theta), \mathcal{I}_c(\theta, A), \mathcal{I}_t(\theta, A)$.

PROPOSITION 2.9. *Under the regularity of the model $\{\mathbb{P}_\theta : \theta \in \Theta\}$ defined by the density functions $\{f_\theta : \theta \in \Theta\}$, the information matrices corresponding to the censored and truncated likelihood functions based on a set $A$ verify the relations*

$$\mathcal{I}_t(\theta, A) = \mathbb{P}_\theta\left(\frac{I_A}{\mathbb{P}_\theta(A)}\left(\frac{\frac{\partial}{\partial\theta}f_\theta}{f_\theta}\right)\left(\frac{\frac{\partial}{\partial\theta}f_\theta}{f_\theta}\right)^T\right) - \frac{(\frac{\partial}{\partial\theta}\mathbb{P}_\theta(A))(\frac{\partial}{\partial\theta}\mathbb{P}_\theta(A))^T}{(\mathbb{P}_\theta(A))^2}, \tag{2.16}$$

$$\begin{aligned}\mathcal{I}_c(\theta, A) &= \mathbb{P}_\theta(A)\mathcal{I}_t(\theta, A) + \frac{(\frac{\partial}{\partial\theta}\mathbb{P}_\theta(A))(\frac{\partial}{\partial\theta}\mathbb{P}_\theta(A))^T}{\mathbb{P}_\theta(A)(1-\mathbb{P}_\theta(A))} \\ &= \mathcal{I}(\theta) - \mathbb{P}_\theta(A^c)\mathcal{I}_t(\theta, A^c).\end{aligned} \tag{2.17}$$

REMARK 2.1. The information matrices above are obtained from different probability models. However, in our setup, censoring or truncation are artificial. This means that, in fact, we will know the size of our data sample, and thus, truncation must be understood as a way of handling the data outside the trimming set, but not as a way of wrongly reconsidering the data size. Therefore, we must take into account the original data size for a correct analysis of the information given from a complete sample through both procedures.

This leads to the consideration of either the *conditional* (to the number of observations that belong to $A$) *truncated information*, or the *expected truncated information*. In the first case, we would associate the information matrix $k_n\mathcal{I}_t(\theta, A)$ to a sample of size $n$ with $k_n$ elements in $A$, while in the second we should associate that given by $n\mathbb{P}_\theta(A)\mathcal{I}_t(\theta, A)$. This last point of view means, in fact, that in our model of complete data the truncated information should be $\mathcal{I}_t^*(\theta, A) = \mathbb{P}_\theta(A)\mathcal{I}_t(\theta, A)$, leading (2.17) to the equivalent relation $\mathcal{I}_c(\theta, A) = \mathcal{I}(\theta) - \mathcal{I}_t^*(\theta, A^c)$. Of course, the Law of Large Numbers guarantees that both definitions give the same asymptotic value.

It should be also stressed that the information obtained with censoring is ever greater than that expected with truncation, as trivially arises from (2.16).

Moreover, (2.16) and Proposition 2.10 also show that in the elliptical model, for sets $A$ taken as (scaled versions of) the MVE of $\mathbb{P}_\theta$, both information matrices coincide for all the parameters related to the location and shape of the distribution, and only differ for the scale parameter. In other words, if we reparameterize $\Sigma$ as $\Xi = \Sigma/|\Sigma|^{1/p}, \varsigma^2 = |\Sigma|^{1/p}$, for the (scaled versions of) the MVE of the elliptical probability $\mathbb{P}_\theta$, the only different component in the information matrices $\mathcal{I}_c(\theta, A)$ and $\mathcal{I}_t^*(\theta, A)$ is that corresponding to the scale parameter $\varsigma$ analyzed in [11].

PROPOSITION 2.10. *Assume regularity of the elliptical model* $\{\mathbb{P}_\theta : \theta \in \Theta\}$. *Let* $\Sigma$ *be reparameterized by* $\Xi = \Sigma/|\Sigma|^{1/p}$, $\varsigma^2 = |\Sigma|^{1/p}$. *Then, for every* $\theta_0 = (\mu_0, \Sigma_0) \in \Theta$ *and every* $r > 0$, *the following relations hold:*

$$\frac{\partial}{\partial \mu}\bigg|_{\theta_0} \mathbb{P}_\theta(\mathcal{E}((\mu, \Sigma), r)) = 0, \qquad \frac{\partial}{\partial \Xi}\bigg|_{\theta_0} \mathbb{P}_\theta(\mathcal{E}((\mu_0, \Sigma_0), r)) = 0.$$

In the GEM, the information matrix $\mathcal{I}_s(\eta, A)$, where $\eta = (\pi, \theta)$, is composed of a sub-matrix corresponding to the parameter $\theta$, a term corresponding to $\pi$ and $p(p+1)/2$ terms (i.e., the same number as the dimension of $\theta$) corresponding to the cross terms between $\theta$ and $\pi$. We will, respectively, denote them by $\mathcal{I}_s(\theta, A)$, $\mathcal{I}_s(\pi, A)$ and $\mathcal{I}_s(\theta^i, \pi, A)$.

PROPOSITION 2.11. *Under the regularity of the model defined by the density functions $\{f_\theta : \theta \in \Theta\}$, the information matrix for the GEM (2.7) verifies*

$$\mathcal{I}_s(\theta, A) = (1-\pi)\left(\mathbb{P}_\theta(A)\mathcal{I}_t(\theta, A) + \frac{(\frac{\partial}{\partial\theta}\mathbb{P}_\theta(A))(\frac{\partial}{\partial\theta}\mathbb{P}_\theta(A))^T}{\mathbb{P}_\theta(A)(1-(1-\pi)\mathbb{P}_\theta(A))}\right),$$

$$\mathcal{I}_s(\pi, A) = \frac{\mathbb{P}_\theta(A)}{(1-\pi)(1-(1-\pi)\mathbb{P}_\theta(A))}, \tag{2.18}$$

$$\mathcal{I}_s(\theta^i, \pi, A) = -\frac{\frac{\partial}{\partial\theta^i}\mathbb{P}_\theta(A)}{1-(1-\pi)\mathbb{P}_\theta(A)}, \qquad i = 1, \dots, p(p+1)/2.$$

REMARK 2.2. Since we will often be interested only in the $\theta$'s parameters, it is natural to explore what is the information for $\theta$ in the GEM (2.7), treating $\pi$ as a nuisance parameter. According to the well-known block matrix form of matrix inverses, the block of the inverse matrix of $\mathcal{I}_s(\eta, A)$ corresponding to the $\theta$'s parameters can be expressed as

$$\left(\mathcal{I}_s(\theta, A) - \frac{\frac{\partial}{\partial\theta}\mathbb{P}_\theta(A)}{1-(1-\pi)\mathbb{P}_\theta(A)}(\mathcal{I}_s(\pi, A))^{-1}\frac{(\frac{\partial}{\partial\theta}\mathbb{P}_\theta(A))^T}{1-(1-\pi)\mathbb{P}_\theta(A)}\right)^{-1}, \tag{2.19}$$

hence, the matrix between the great parentheses is considered as the information for $\theta$.

From (2.18), it is straightforward that this information coincides with $(1-\pi)\mathbb{P}_\theta(A)\mathcal{I}_t(\theta, A)$. This agrees with our considerations in Remark 2.1 and the second item in Proposition 2.7: Taking into account that the truncated model associated to the GEM on $A$ coincides with the truncated model (on $A$) associated to the uncontaminated model $\{\mathbb{P}_\theta : \theta \in \Theta\}$, the information matrix for the GEM must coincide with that of the truncated model corrected through a suitable factor. The expected number of sample data points in $A$ obtained from a random sample of size $n$ from the contaminated probability $\mathbb{P}^*_{\theta_0} = (1-\pi_0)\mathbb{P}_{\theta_0} + \pi_0 Q$, where $Q$ is any probability with support in $A^c$, is precisely $n(1-\pi_0)\mathbb{P}_{\theta_0}(A)$, thus, the truncated information obtained from one observation from the original GEM should be

$$\mathcal{I}_t^{**}(\theta_0, A) = (1-\pi_0)\mathbb{P}_{\theta_0}(A)\mathcal{I}_t(\theta_0, A). \tag{2.20}$$

This also supports that the MLE(s) coincides with the MLE(t) asymptotically.

**3. Robustness and asymptotics of the estimators.** In the finite sample setting, the robustness of an estimator is usually measured through its (finite sample) BP, which for an estimator $T_n$ based on a sample $\mathcal{X}_n$ will be denoted as $\varepsilon^*(T_n, \mathcal{X}_n)$. Of course, the BP has no sense if we are only able to assure the asymptotic existence of an estimator. In fact, its analysis is closely related to arguments on the existence of the estimator. In our case we have shown in Propositions 2.5 and 2.6 the existence of the MLE(s), MLE(c) and MLE(r) under very general hypotheses. If our initial estimator is equivariant, the ellipsoid on which we base our ML (final) estimation will be also equivariant and the whole procedure will obviously maintain the equivariance property. But, as stated in Theorem 3.1, our one-step procedures also preserve the initial BP. In fact, by merging the arguments in Section 5 in [19] with those used in the discussion showing the existence of our estimators, it is straightforward to show the following theorem.

THEOREM 3.1. *Let $\mathcal{X} = \{x_1, \dots, x_n\} \subset \mathbb{R}^p$, $n > p$, be a sample of points in general position. Let $t_n$ and $C_n$ be estimates of location and covariance. Let*

$$A_n := \{x \in \mathbb{R}^p : (x - t_n)^T C_n^{-1} (x - t_n) \le c_1\},$$

*where $c_1$ is any fixed value such that the set $A_n$ contains at least $[\frac{n+p+1}{2}]$ points of $\mathcal{X}$.*

*If the hypotheses in Proposition 2.5 are satisfied and $\hat{\theta}_{s,n}$, $\hat{\theta}_{c,n}$ and $\hat{\theta}_{r,n}$ are respectively the MLE(s), MLE(c) and MLE(r) based on $A_n$, then*

$$\min\{\varepsilon^*(\hat{\theta}_{s,n}, \mathcal{X}), \varepsilon^*(\hat{\theta}_{c,n}, \mathcal{X}), \varepsilon^*(\hat{\theta}_{r,n}, \mathcal{X})\} \ge \min\{\varepsilon^*(t_n, \mathcal{X}), \varepsilon^*(C_n, \mathcal{X})\}.$$

*In particular, when $t_n$ and $C_n$ are the MVE-based estimators, then*

$$\varepsilon^*(\hat{\theta}_{s,n}, \mathcal{X}) = \varepsilon^*(\hat{\theta}_{c,n}, \mathcal{X}) = \varepsilon^*(\hat{\theta}_{r,n}, \mathcal{X}) = [(n - p + 1)/2].$$

In order to obtain the Influence Functions (IF) of our estimators, we will begin with a fixed ellipsoid $A = \mathcal{E}(\gamma)$, $\gamma \in \Gamma$, and emphasize on the dependence on the parameter $\gamma$. In this case the IF's of our estimators can be obtained as the IF's of M-estimators. Thus, after Section 2.3, under the usual conditions to allow for interchanging differentiation and integration [recall relation (2.15) obtained from the regularity of the model], and under the assumed model, provided that the involved information matrices are nonsingular, we obtain

$$\mathrm{IF}(x, \hat{\theta}_{*,n}(\gamma), \theta_0) = -\left(\mathbb{P}_{\theta_0}\left(\frac{\partial}{\partial\theta} h^*_{\theta,\gamma}\right)\right)^{-1} h^*_{\theta_0,\gamma}(x), \tag{3.1}$$

where $\hat{\theta}_{*,n}(\gamma) = \hat{\theta}_{t,n}(\gamma)$ or $\hat{\theta}_{c,n}(\gamma)$ and $h^*_{\theta,\gamma} = h^t_{\theta,\gamma}$ or $h^c_{\theta,\gamma}$, defined by

$$h^t_{\theta,\gamma} := \left(\frac{\frac{\partial}{\partial\theta} f_\theta}{f_\theta} - \frac{\frac{\partial}{\partial\theta}\mathbb{P}_\theta(\mathcal{E}(\gamma))}{\mathbb{P}_\theta(\mathcal{E}(\gamma))}\right) I_{\mathcal{E}(\gamma)},$$

$$h^c_{\theta,\gamma} := \frac{\frac{\partial}{\partial\theta} f_\theta}{f_\theta} I_{\mathcal{E}(\gamma)} - \frac{\frac{\partial}{\partial\theta}\mathbb{P}_\theta(\mathcal{E}(\gamma))}{\mathbb{P}_\theta(\mathcal{E}(\gamma)^c)} I_{\mathcal{E}(\gamma)^c}.$$

On the other hand, under the GEM, by defining

$$h^s_{\theta,\pi,\gamma} := \left( \frac{\frac{\partial}{\partial\theta} f_\theta}{f_\theta} I_{\mathcal{E}(\gamma)} + \frac{(1-\pi)\frac{\partial}{\partial\theta}\mathbb{P}_\theta(\mathcal{E}(\gamma)^c)}{(1-\pi)\mathbb{P}_\theta(\mathcal{E}(\gamma)^c)+\pi} I_{\mathcal{E}(\gamma)^c}, \right.$$
$$\left. \frac{-1}{1-\pi} I_{\mathcal{E}(\gamma)} + \frac{1-\mathbb{P}_\theta(\mathcal{E}(\gamma)^c)}{1-(1-\pi)\mathbb{P}_\theta(\mathcal{E}(\gamma))} I_{\mathcal{E}(\gamma)^c} \right)^T,$$

and recalling the information matrix $\mathcal{I}_s(\eta, \mathcal{E}(\gamma))$ (see Proposition 2.11), we obtain

$$\mathrm{IF}(x, \hat{\theta}_{s,n}(\gamma), \hat{\pi}_{s,n}(\gamma), \theta_0, \pi_0) = (\mathcal{I}_s(\eta_0, \mathcal{E}(\gamma)))^{-1} h^s_{\theta_0,\pi_0,\gamma}(x). \tag{3.2}$$

Because of the continuity of the estimators with respect to $\gamma$, it is easy to see that the IF of the estimator $\hat{\theta}_{*,n}(\gamma_n)$ coincides with that of $\hat{\theta}_{*,n}(\gamma)$ if $\{\gamma_n\}_n \subset \Gamma$ and $\gamma_n \to \gamma \in \Gamma$, if we apply the main idea in the proof of Theorem B.1 in [10] to the points that do not belong to the boundary of $\mathcal{E}(\gamma)$. Therefore, the IF of the one-step (truncated, censored or smart) estimator based on the MVE estimators will be the one given by (3.1), or (3.2) with $\mathcal{E}(\gamma)$ being the MVE of $\mathbb{P}$, where $\mathbb{P}$ belongs to an elliptical model [or to the GEM model given by (2.2) for the elliptical model].

Of course, the asymptotic variances computed from the information matrices $\mathcal{I}_c(\theta_0, \mathcal{E}(\gamma))$, $\mathcal{I}^*_t(\theta_0, \mathcal{E}(\gamma))$ and $\mathcal{I}^{**}_t(\theta_0, \mathcal{E}(\gamma))$, taking into account Remarks 2.1 and 2.2, and by integration of the square of the relations (3.1) and (3.2), coincide.

3.1. *Strong consistency.* To explore the asymptotic behavior of our estimators, we begin with the consistency of the initial estimator. We will show that any initial consistent estimator under the model would give the same asymptotic behavior. The consistency of the MVE in the uncontaminated model has already been treated in [7]. However, under the uniqueness of the theoretical MVE, it is not difficult to show the following proposition that covers the GEM.

PROPOSITION 3.1. *Let $g$ be a decreasing function which defines an elliptical family on $\mathbb{R}^p$. Let $\{X_n\}_n$ be a random sample obtained from the distribution*

$$\mathbb{P} = (1-\pi_0)\mathbb{P}_{\theta_0} + \pi_0 Q \tag{3.3}$$

*in the GEM of the elliptical family defined by $g$ with $\pi_0 < 1/2$, let $A = \mathcal{E}(\mu, \Sigma, r)$ be the MVE of $\mathbb{P}$, which we assume to be unique, and $A_n = \mathcal{E}(\mu_n, \Sigma_n, r_n)$ be the sample MVE. Then we have that $\lim_n I_{A_n} = I_A$ a.s.*

Now, we are in a position to prove the consistency of the smart estimate under the GEM.

THEOREM 3.2 (Consistency of the estimators). *Let $g$ be a function which defines an elliptical family on $\mathbb{R}^p$ and satisfies* G1, G2 *and* Gp. *Let $\{X_n\}_n$ be a random sample taken from the distribution*

$$\mathbb{P} = (1-\pi_0)\mathbb{P}_{\theta_0} + \pi_0 Q$$

*in the GEM of the elliptical family defined by $g$, $0 \le \pi_0 < 1/2$ and $\theta_0 \in \Theta$.*

*Let $A$ be an MVE of $\mathbb{P}$, which we assume to be unique. Let $\{A_n\}_n$ be a sequence of empirical MVE's and let $\{(\hat{\pi}_n, \hat{\mu}_n, \hat{\Sigma}_n)\}_n$ be a sequence of MLE(s) based on the ellipsoids $\{A_n\}_n$. Then, the following is satisfied:*

1. *The MLE(s) based on $\{A_n\}_n$ is strongly consistent.*
2. *If $\pi_0 = 0$, then the MLE(c) based on $\{A_n\}_n$ is strongly consistent.*
3. *If $\alpha \in (0, 1/2)$, then the MLE(r), $\hat{\theta}_{r,n}$, based on $\{A_n\}_n$ under the restrictions given by $\Theta^{\alpha}_{A_n}$ is strongly consistent. Moreover, if $\hat{\pi}_n$ is computed from $\hat{\theta}_{r,n}$ using (2.14), then also $\hat{\pi}_n \to \pi_0$ a.s.*

Proposition 3.2 shows that, for a large enough sample size, the restricted parameter set contains the true value of the parameter. Thus, for large sizes these restrictions are, in fact, superfluous.

PROPOSITION 3.2. *Assume the hypotheses of Theorem 3.2. Let $\alpha \in (0, 1/2)$ be given and let $\Theta^{\alpha}_{A_n}$ be defined as in (2.13). Then, for a.e. sample there exists $\delta > 0$ such that $\{\theta : \|\theta - \theta_0\| < \delta\} \subset \Theta^{\alpha}_{A_n}$, for large enough $n$.*

3.2. *Asymptotic distribution.* Although the extension of the argmax-based arguments of the Empirical Processes Theory to the semiparametric framework is certainly not trivial, our model is well suited for such a task, because of the special features of the family of ellipsoids parameterized through the set $\Gamma$. In fact, Section 3.2.4 of [26] can easily be tuned to cover our setup by verbatim repeating the reasoning therein in order to get the chain of results on linearization given in the Appendix as well as their consequences.

We only consider with some detail the MLE(r) which needs some additional analysis. Let $\alpha \in (0, 1/2)$, and $\{\gamma_n\}$ be the sequence of parameters associated to the sequence of sample MVE's. We initially assume the hypotheses in Lemma A.7, as well as the regularity of the underlying elliptical model. After the consistency results, for the analysis of the asymptotic distribution, we can assume that the $\gamma$-parameters belong to a compact subset $K$ of $\Gamma$, and that the $\theta$-parameters verify the restrictions given by $\Theta^{\alpha}_{\mathcal{E}(\gamma_n)}$ and belong to the set $\{\theta : \|\theta - \theta_0\| < \delta\}$ for some $\delta > 0$ and large enough $n$.

Let us consider the function $m_{\theta,\gamma}$, associated to the MLE(r), given by

$$m_{\theta,\gamma}(x) := I_{\mathcal{E}(\gamma)} \log\left(\frac{g((x-\mu)^T \Sigma^{-1}(x-\mu))}{\int_{\mathcal{E}(\gamma)} g((y-\mu)^T \Sigma^{-1}(y-\mu))\, dy}\right). \tag{3.4}$$

Lemma 3.1 allows us to apply Theorem A.1 under the condition required in Lemma A.8.

LEMMA 3.1. *Let us assume that $g$ is twice continuously differentiable. Let $\theta_0 = (\mu_0, \Sigma_0) \in \Theta$, $\gamma_0 \in \Gamma$ be such that*

$$\inf_{x \in \mathcal{E}(\gamma_0)} g((x-\mu_0)^T \Sigma_0^{-1}(x-\mu_0)) > 0.$$

*Then, there exist a vector valued function $\dot{m}_{\theta\gamma}$, $\delta > 0$ and a compact neighborhood $K$ of $\gamma_0$ such that*

$$\left\{ \frac{m_{\theta\gamma} - m_{\theta_0\gamma} - (\theta-\theta_0)^T \dot{m}_{\theta_0\gamma}}{\|\theta-\theta_0\|} : \|\theta-\theta_0\| < \delta, \gamma \in K \right\} \tag{3.5}$$

*is $\mathbb{P}$-Donsker and*

$$\mathbb{P}(m_{\theta\gamma} - m_{\theta_0\gamma} - (\theta-\theta_0)^T \dot{m}_{\theta_0\gamma})^2 = o(\|\theta-\theta_0\|)^2, \tag{3.6}$$

*uniformly in $\gamma \in K$.*

If the matrix of second derivatives is continuous and nonsingular, relation (A.24) in Lemma A.8 and the consistency of the MVE's produce the asymptotic laws of the MLE(s) at the announced rate $n^{1/2}$, independently of the rate of convergence of the initial estimator.

The analogous result for the MLE(c) under the elliptical model follows from similar considerations. Finally, under a probability in the GEM of the elliptical model with $\pi > 0$, the consistency of the MLE(s) assures, from Proposition 2.7, that the MLE(s) coincides with the natural MLE(r) asymptotically. Thus, they share their asymptotic normal distribution, with the covariance matrices related to the information matrices already obtained. We summarize these results in the following final theorem.

THEOREM 3.3 (Asymptotic distributions). *Assume the hypothesis in Lemma A.7 and that $g$ is twice continuously differentiable. Let $A$ be the (only) MVE of $\mathbb{P}$. For each $n \in \mathbb{N}$ , consider an estimation, $A_n$, obtained through a consistent estimator of the MVE; the MLE(r), $\hat{\theta}_{r,n}$ under the restriction defined by $\Theta^{\alpha}_{A_n}$ for some $\alpha \in (0, 1/2)$, as well as the MLE(s), $\hat{\theta}_{s,n}$, and the MLE(c), $\hat{\theta}_{c,n}$, based on $A_n$.*

*If the corresponding information matrices are nonsingular, then:*

1. *$\sqrt{n}(\hat{\theta}_{r,n} - \theta_0)$ converges in law to a centered multivariate normal distribution with covariance matrix given by the inverse of the information matrix, $\mathcal{I}_t^{**}(\theta, A)$ defined through (2.20) and Proposition 2.9.*
2. *If $\pi_0 = 0$, $\sqrt{n}(\hat{\theta}_{c,n} - \theta_0)$ converges in law to a centered multivariate normal distribution with covariance matrix given by the inverse of the information matrix, $\mathcal{I}_c(\theta, A)$, defined in Proposition 2.9.*

3. *If* $\pi_0 > 0$, $\sqrt{n}(\hat{\theta}_{s,n} - \theta_0)$ *converges in law to a centered multivariate normal distribution with covariance matrix given in (2.19).*

**4. Discussion.** A consideration on the efficiency of the obtained estimators can be illuminating. Note that the rate of convergence is always $n^{1/2}$, but also that the asymptotic law of the estimators depends on the limit ellipsoid but not on the rate of convergence of the initial estimator to this ellipsoid. In fact, from the asymptotic results and the expressions of the information matrices, it becomes apparent that the efficiency is equivalent to that obtained from the corresponding MLE computed on the theoretical (enlarged) MVE. Therefore, it is greater than that obtained by the usual one step reweighting, even for initial estimators that converge faster than the MVE estimator.

Under the elliptical model, any high-BP consistent initial estimator $(t_n, C_n)$ of $\theta = (\mu, \Sigma)$ could be used to produce our estimation $A_n := \mathcal{E}(t_n, C_n, r_{C_n,\alpha})$ of the central ellipsoid $A = \mathcal{E}(\mu, \Sigma, r_{\Sigma,\alpha})$, covering, say, the $1 - \alpha = 95\%$ of the theoretical distribution. Between our proposals, the MLE(c) based on $A_n$ will provide maximum efficiency and the same BP as $(t_n, C_n)$. We recall that, according to Remark 2.1, the gain of efficiency with respect to the MLE(t) appears only in the estimation of $\Sigma$, thus, the MLE(t) and MLE(c) of $\mu$ based on $A_n$ have the same efficiency. Since it is usual to justify the use of robust estimators looking at the behavior under the model (i.e., assuming the existence of no contamination), the greater efficiency of the MLE(c) under the elliptical model would justify its prioritary use.

In Tables 1, 2 and 3 (see Appendix B) we present the asymptotic efficiencies of these estimators under the uncontaminated elliptical model and their versions based on enlarged MVE estimations. The comparison of the efficiencies in these tables with other well-known robust estimators (see, for instance, the efficiencies obtained in [5]) or [4] shows that the combination "Initial MVE estimator"+"Scaled version for a given $\alpha$" + "MLE(c)" gives better efficiencies between the highest-BP equivariant estimators.

In the contaminated model, it is intuitively sound that the best choice for an estimator based on a subsample which contains no outlier should be the MLE(t). In this sense the MLE(s) is the natural MLE, because it only substitutes the estimation provided by the MLE(t) when it does not exist or the sample does not sufficiently match the GEM. In some way it also robustifies the MLE(t) that, as already noted, possibly does not even exist. This nonexistence is in apparent contradiction with Theorem 1 in [20], but the BP studied in this theorem is not the sample-based one and does not reflect the possible nonexistence of the MLE(t), which could make it undesirable from the robustness point of view. The MLE(r) would be an excellent alternative, taking into account the choice of the initial trimmed sets.

In the presence of outliers, our choice of the MVE as initial estimator to produce the trimmed set is related to the GEM model, which is based on the possibility of discarding the outliers by resorting to a common *central ellipsoid* of the contaminated and uncontaminated models. Since our proposals circumvent the drawback of its convergence rate, this choice stresses the improvement of efficiency obtained through the presented methodology.

The literature on robust estimation in the elliptical model usually analyzes the estimation of $|\Sigma|$ *a posteriori*, by adjusting on the basis of the model and the estimates of location and shape. This generally leads to a Fisher inconsistent estimation under a real contaminated model, even in the considered GEM in which only outliers contaminate the distribution. On the contrary, our proposals are in their own right MLE, even for the size of $\Sigma$, and only the MLE(c) would be Fisher inconsistent (when $\pi_0 > 0$).

In the applications, every proposal can be computed through a variant of the EM algorithm (see Section 4.2 in Dempster, Laird and Rubin [9]) and based on the improved MVE given in (6.59) in [22]. The variant of the EM algorithm can be based on a Monte Carlo approximation to the integrals using a random sample from the appropriate elliptical distribution, while in the M step we need to solve the estimation problem for the original (nontruncated, noncensored) elliptical distribution.

## APPENDIX A: PROOFS AND SOME TECHNICAL RESULTS

Proof of Proposition 2.1. If $f_{\theta_1} = f_{\theta_2}$ $\lambda$-a.s. on $A$, then

$$\lambda\left\{x \in A : \sum_{i=1}^{j}(Q_j(\theta_1) - Q_j(\theta_2))T_j(x) = \log(C(\theta_1)/C(\theta_2))\right\} = \lambda(A),$$

and the $T$'s would satisfy a linear constraint on $A$ with $\lambda$-positive measure. □

Proof of Proposition 2.2. Let $\theta = (\mu, \Sigma) \in \Theta$ be such that, for some $k$, it satisfies

$$\mathbb{P}_{\theta_0}[C_k] = \mathbb{P}_{\theta_0}[A], \tag{A.1}$$

where $C_k = \{x \in A : f_\theta(x) = k f_{\theta_0}(x)\}$.

Assume that $\mu \neq \mu_0$ and let $\varepsilon > 0$ such that $B(\mu_0, \varepsilon) \subset A$ and $f_{\theta_0} > 0$ on $B(\mu_0, \varepsilon)$. Then $\varepsilon^* = \inf(\|\mu - \mu_0\|, \varepsilon) > 0$. For every $x \in S(\mu_0, \varepsilon^*) \cap \{y : \langle y - \mu, \mu_0 - \mu\rangle > 0\}$, let $[x, \mu_0] \subset \mathbb{R}^p$ be the segment joining $x$ with $\mu_0$.

The function $f_{\theta_0}$ (resp. $f_\theta$) increases (resp. decreases) on the segment from $x$ to $\mu_0$. Because of G1, $f_{\theta_0}$ is not constant on this segment. Thus, if we denote $\lambda_x^1$ the (one-dimensional) Lebesgue measure on $[x, \mu_0]$, then

$$\lambda_x^1\{f_\theta \neq k f_{\theta_0}\} \cap [x, \mu_0] > 0,$$

which makes (A.1) impossible.

This implies that $\mu = \mu_0$. Moreover, from (A.1) there exists a sequence $\{x_n\} \subset C_k - \{0\}$, such that $\lim_n x_n = \mu_0$ and, since (from G1) $g(0+) > 0$, we obtain

$$k = \lim_n \frac{f_\theta(x_n)}{f_{\theta_0}(x_n)} = \left(\frac{|\Sigma_0|}{|\Sigma|}\right)^{1/2} \lim_n \frac{g((x_n-\mu_0)^T\Sigma^{-1}(x_n-\mu_0))}{g((x_n-\mu_0)^T\Sigma_0^{-1}(x_n-\mu_0))} = \left(\frac{|\Sigma_0|}{|\Sigma|}\right)^{1/2}.$$

In other words, if $x \in C_k$, we have

$$g((x-\mu_0)^T\Sigma^{-1}(x-\mu_0)) = g((x-\mu_0)^T\Sigma_0^{-1}(x-\mu_0)).$$

On the other hand, let $x \in A$ such that $f_{\theta_0}(x) < f_{\theta_0}(\mu_0)$ and let

$$t_x = \sup\{t > 0 : g(t) > |\Sigma_0|^{1/2} f_{\theta_0}(x)\}.$$

Because of G1, $t_x \neq 0$. Moreover, taking into account that $A$ is open and (A.1), we have that there exist two sequences $\{x_n\}, \{y_n\} \subset C_k$ such that

$$\lim_n (x_n-\mu_0)^T\Sigma_0^{-1}(x_n-\mu_0) = \lim_n (y_n-\mu_0)^T\Sigma_0^{-1}(y_n-\mu_0) = t_x, \tag{A.2}$$

while $(y_n-\mu_0)^T\Sigma_0^{-1}(y_n-\mu_0) > t_x$, and $(x_n-\mu_0)^T\Sigma_0^{-1}(x_n-\mu_0) < t_x$, for every $n \in \mathbb{N}$. Therefore, by definition of $C_k$ and $t_x$, it must also happen that

$$\lim_n (y_n-\mu_0)^T\Sigma^{-1}(y_n-\mu_0) = t_x. \tag{A.3}$$

Without loss of generality, we can assume that the sequence $\{y_n\}_n$ is convergent. Let $y_x$ be its limit. Thus, from (A.2) and (A.3) we have that

$$(y_x-\mu_0)^T\Sigma^{-1}(y_x-\mu_0) = t_x \quad \text{and} \quad (y_x-\mu_0)^T\Sigma_0^{-1}(y_x-\mu_0) = t_x. \tag{A.4}$$

However, by G1, it is possible to choose $x$ in order to obtain infinite different values for $t_x$ above. This and the freedom we have to choose the convergent sequence $\{y_n\}$ give that at most there exists a matrix $\Sigma$ which satisfies the infinite number of relations included in (A.4). Since $\Sigma_0$ satisfies all these relations, the only possibility is to have $\Sigma = \Sigma_0$. $\square$

Proposition 2.4 employs the following lemma in its proof.

LEMMA A.1. *Let $\mathcal{X} := \{x_1, \dots, x_n\}$, where $n > p$, be a set whose points are in general position. Let $\mathcal{H}$ be the family of all hyperplanes in $\mathbb{R}^p$, and given $H \in \mathcal{H}$, let us denote the distance from $x_i$ to $H$ by $d_i(H) := \inf\{\|x_i - h\| : h \in H\}$, $i = 1, \dots, n$.*

*If $(d_{(1)}(H), d_{(2)}(H), \dots, d_{(n)}(H))$ is the ordered set of the values, $d_i(H)$, $i = 1, \dots, n$, then $\inf_{H\in\mathcal{H}} d_{(p+1)}(H) > 0$.*

PROOF. Every $H \in \mathcal{H}$ is determined by the vector $v \in S_{p-1}$, the unit sphere in $\mathbb{R}^p$, and the value $b \in \mathbb{R}$ which satisfy $H = \{x \in \mathbb{R}^p : \langle x, v \rangle = b\}$. Let us denote $H = H_{v,b}$. Also, set $x_i = (x_i^1, \ldots, x_i^p)$, for every $i = 1, \ldots, n$, $M_n := \sup\{|x_i^j| : j = 1, \ldots, p \text{ and } i = 1, \ldots, n\}$ and $\mathcal{H}_n = \{(v, b) \in S_{p-1} \times \mathbb{R} : H_{v,b} \cap [-M_n, M_n]^p \neq \varnothing\}$. Obviously,

$$\inf\{d_{(p+1)}(H) : H \in \mathcal{H}\} = \inf\{d_{(p+1)}(H_{v,b}) : (v, b) \in \mathcal{H}_n\}.$$

For every $v \in S_{p-1}$ such that there exists $b \in \mathbb{R}$, which satisfies that $(v, b) \in \mathcal{H}_n$, let us consider the continuous maps

$$v \to B^n(v) := \sup\{b \in \mathbb{R}^p : (v, b) \in \mathcal{H}_n\},$$
$$v \to B_n(v) := \inf\{b \in \mathbb{R}^p : (v, b) \in \mathcal{H}_n\}.$$

Since $S_{p-1}$ is compact and $\mathcal{H}_n = \bigcup_{v \in S_{p-1}} \{v\} \times [B_n(v), B^n(v)]$, where we take $[B^n(v), B_n(v)] = \varnothing$ if these maps are not defined, we obtain that $\mathcal{H}_n$ is compact.

On the other hand, for every $i = 1, \ldots, n$, the map $(v, b) :\to d_i(H_{v,b})$ is continuous, hence, $(v, b) :\to d_{(p+1)}(H_{v,b})$ is also continuous and reaches its infimum on $\mathcal{H}_n$, proving the lemma from the general position assumption. □

PROOF OF PROPOSITION 2.4. We will only consider the more involved case $p > 1$. Let $\{(\mu_k, \Sigma_k)\}_k \subset \Theta$ be a sequence such that

$$\lim_k \prod_{i=1}^n f_{(\mu_k, \Sigma_k)}(x_i) = \sup_{(\mu, \Sigma) \in \Theta} \prod_{i=1}^n f_{(\mu, \Sigma)}(x_i). \tag{A.5}$$

Since $g$ is continuous and $g(0) > 0$, it must be

$$\liminf_k \prod_{i=1}^n f_{(\mu_k, \Sigma_k)}(x_i) > 0. \tag{A.6}$$

Let $v_k^1, \ldots, v_k^p$ and $\delta_k^1, \ldots, \delta_k^p$ be the eigenvectors and eigenvalues of $\Sigma_k$. Let $\Delta_k := \inf\{\delta_k^1, \ldots, \delta_k^p\}$ and let $j_k$ be such that $\Delta_k = \delta_k^{j_k}$. First, we prove that it is impossible that $\liminf_k \Delta_k = 0$. Let us assume that, on the contrary, there exists a subsequence, which we will denote as the original one, such that $\lim_k \Delta_k = 0$.

Since the points in $\mathcal{X}$ are in general position and, since $n > p$, we can apply Lemma A.1 to obtain that there exists $d > 0$ such that

$$I_k := \{i \in \{1, \ldots, n\} : |(x_i - \mu_i)^T v_i^{j_k}| \geq d\}$$

is a set which contains at least $(n - p)$ elements. Therefore, if $i \in I_k$, then

$$(x_i - \mu_k)^T \Sigma_k^{-1} (x_i - \mu_k) \geq |(x_i - \mu_i)^T v_i^{j_k}|^2 (\Delta_k)^{-1} \geq d^2 (\Delta_k)^{-1},$$

and, since $g$ is nonincreasing, if $i \in I_k$, and $|\Sigma_k| \geq (\Delta_k)^p$, we have that

$$\prod_{i=1}^{n} f_{(\mu_k,\Sigma_k)}(x_i) \leq (\Delta_k)^{-pn/2} g(0)^p g(d^2(\Delta_k)^{-1})^{n-p}. \tag{A.7}$$

Thus, applying assumption Gp in (A.7), we have that if $k$ is big enough,

$$\prod_{i=1}^{n} f_{(\mu_k,\Sigma_k)}(x_i) \leq (\Delta_k)^{n(2\gamma-p)/2-p\gamma} g(0)^p d^{-2\gamma(n-p)},$$

which converges to zero as $n \to \infty$ and contradicts (A.6).

Now, let $\Delta^k := \sup\{\delta_k^1, \dots, \delta_k^p\}$. Since $|\Sigma_k| \geq \Delta_k^{p-1}\Delta^k$, we have that

$$\prod_{i=1}^{n} f_{(\mu_k,\Sigma_k)}(x_i) \leq (\Delta_k)^{-n(p-1)/2}(\Delta^k)^{-n/2} g(0)^n,$$

and, to avoid contradictions with (A.6), we obtain that $\limsup_k |\Delta^k| < \infty$.

Because $\liminf_k |\Delta_k| > 0$ and $\limsup_k |\Delta^k| < \infty$, we can conclude that $\limsup_k \|\mu_k\| < \infty$ because, on the contrary, (A.6) would be false.

This means that every sequence which satisfies (A.5) is contained in a compact set and, in consequence, it contains a convergent subsequence to, say, $(\hat{\mu}, \hat{\Sigma}) \in \Theta$. An easy argument of continuity shows now that $(\hat{\mu}, \hat{\Sigma})$ is the point we are looking for. □

PROOF OF PROPOSITION 2.5. Let $\{(\theta_k, \pi_k)\}_k$ in $\Theta \times [0,1)$, $\theta_k = (\mu_k, \Sigma_k)$, be a sequence such that

$$\lim_{k\to\infty} P_n L^s_{\theta_k,\pi_k/A} = \sup_{\theta\in\Theta,\pi\in[0,1)} P_n L^s_{\theta,\pi/A}.$$

Since $[0,1)$ is bounded, we can assume that there exists $\hat{\pi} = \lim_k \pi_k$.

Taking into account that the second summand in (2.8) is bounded above, we can repeat the same reasoning as in Proposition 2.4 to show that there exists a convergent subsequence of $\{\theta_k\}$ whose limit belongs to $\Theta$ and also that $\hat{\pi} < 1$. Thus, from the continuity of $f_\theta$ and $\mathbb{P}_\theta$ [recall (2.12)], we obtain that the maximum is attained at the limit of this subsequence.

The proof for the MLE(c) is the same, by keeping $\pi = 0$ fixed. □

PROOF OF PROPOSITION 2.6. This proof goes along the same lines as the one we gave for Proposition 2.4 because, under the restrictions, the term $\mathbb{P}_{\theta_k}(A)$ is bounded away from zero. □

PROOF OF PROPOSITION 2.7. The first statement directly follows from the expression (2.8) of the objective function, which, for $\pi = 0$, coincides with that of the censored framework. Concerning the other statement, notice that an equivalent expression for (2.8) is

$$P_n I_A \log(f_\theta/\mathbb{P}_\theta(A)) + P_n(I_A \log((1-\pi)\mathbb{P}_\theta(A)) + I_{A^c}\log(1-(1-\pi)\mathbb{P}_\theta(A))). \tag{A.8}$$

Let us denote $\psi(\pi,\theta)$ to the sum of the second and third summands in this expression (which are the only ones depending on $\pi$). Note that $\psi(\pi^*(\theta),\theta)$ does not depend on $\theta$. On the other hand, derivation of $\psi$ with respect to $\pi$ easily shows that, for every $\theta$, if $\pi \geq 0$ and $\pi > \pi^*(\theta)$, then $\psi(\pi,\theta)$ is nonincreasing on $\pi$, thus, the maximum value of $\psi(\pi,\theta)$ on $[0,1]$ is $\psi(\pi^*(\theta),\theta)$ if $\pi^*(\theta) \geq 0$ else $\psi(0,\theta)$. Then it follows that $\hat{\pi}_{s,n} = \pi^*(\hat{\theta}_{s,n})$ when $\hat{\pi}_{s,n} > 0$, and the maximum value of (A.8) under the restriction $\mathbb{P}_\theta(A) \geq P_n(A)$ is, as stated in the second item,

$$P_n I_A \log(f_\theta/\mathbb{P}_\theta(A)) + \psi(\pi^*(\theta),\theta). \qquad \square$$

The proof of Theorem 2.1 is based on Lemma A.2. Let $\mu$ be a positive $\sigma$-finite measure on $\beta^p$ such that the function $c$ on $\mathbb{R}^p$ defined by $c(\theta) = \int \exp\{\sum_{i=1}^d \theta_j x_j\}\mu(dx)$ is not identically $+\infty$. $c$ is the so-called *Laplace transform of* $\mu$. Its domain is the set $\{\theta \in \mathbb{R}^d : c(\theta) < +\infty\}$.

LEMMA A.2 (Theorem 7.1 in [1]). *Let $\kappa = \log c$ be the logarithm of the Laplace transform of $\mu$. Then $\kappa$ is a closed convex function on $\mathbb{R}^p$ and is strictly convex on its domain, provided $\mu$ is not concentrated on an affine subspace of $\mathbb{R}^p$.*

PROOF OF THEOREM 2.1. We will employ the canonical form of (2.3), obtained by a re-parameterization and the absorption of $h$ into $\lambda$, leading to

$$f_\theta(x) = C(\theta)\exp\left\{\sum_{j=1}^d \theta_j T_j(x)\right\}. \tag{A.9}$$

The expression within the brackets in (2.6) is the logarithm of a density function, say, $g_{\theta,\pi}$, with respect to the $\sigma$-finite measure $\lambda|_A + \delta_{\{\mathbf{c}\}}$,

$$g_{\theta,\pi}(x) := I_A(x)(1-\pi)f_\theta(x) + I_{\{\mathbf{c}\}}(x)[(1-\pi)(1-\mathbb{P}_\theta(A)) + \pi]. \tag{A.10}$$

It is straightforward to obtain the following exponential expression for $g_{\theta,\pi}$, whenever the condition $1-(1-\pi)\mathbb{P}_\theta(A) > 0$ holds [or in an equivalent way, whenever $-\pi < \mathbb{P}_\theta(A^c)/\mathbb{P}_\theta(A)$, allowing even negative values of $\pi$]:

$$g_{\theta,\pi} = (1-(1-\pi)\mathbb{P}_\theta(A))\left[I_A \exp\left\{\sum_{j=1}^d \theta_j T_j\right\}\frac{(1-\pi)C(\theta)}{1-(1-\pi)\mathbb{P}_\theta(A)} + I_{\{\mathbf{c}\}}\right]$$

$$= (1-(1-\pi)\mathbb{P}_\theta(A)) \exp\left\{ I_A \left( \sum_{j=1}^d \theta_j T_j + \log \frac{(1-\pi)C(\theta)}{1-(1-\pi)\mathbb{P}_\theta(A)} \right) \right\},$$

easily seen within the class of exponential distributions, if we add the parameter

$$\theta_{d+1} = \log\left( \frac{(1-\pi)C(\theta)}{(1-(1-\pi)\int_A C(\theta)\exp\{\sum_{j=1}^d \theta_j T_j(x)\}\lambda(dx))} \right), \tag{A.11}$$

leading to the density functions with respect to $\lambda|_A + \delta_{\{\mathbf{c}\}}$,

$$h_{\theta,\theta_{d+1}}(x) := D(\theta, \theta_{d+1}) \exp\left\{ \sum_{j=1}^d \theta_j (T_j(x) I_A(x)) + \theta_{d+1} I_A(x) \right\}. \tag{A.12}$$

The hypothesis requiring that the $T$'s do not satisfy a linear constraint on $A$ implies that $T_1, T_2, \dots, T_d, I_A$ also do not satisfy such a linear constraint on $A$. This allows us, by Lemma A.2, to guarantee that $-\log D(\theta, \theta_{d+1})$ is a closed strictly concave function on its domain.

On the other hand, from (A.11), the restriction $\pi \geq 0$ can be written as

$$\theta_{d+1} + \log\left( \int \exp\left\{ \sum_{j=1}^d \theta_j T_j(x) \right\} \lambda^*(dx) \right) \leq 0, \tag{A.13}$$

where the term $I_{A^c}$ is included in a new measure $\lambda^* = \lambda|_{A^c}$.

Once more by Lemma A.2, the function on the left of (A.13) is a convex function. Therefore, the restricted set defined by (A.13) is a convex set.

Let $\mathbb{P}$ be a probability distribution verifying the hypotheses. The function

$$\mathbb{P} \log h_{\theta,\theta_{d+1}} = \sum_{j=1}^d \theta_j \mathbb{P}(T_j I_A) + \theta_{d+1} \mathbb{P}(I_A) - \log D(\theta, \theta_{d+1})$$

is then a strictly concave function on its domain, so, if any, it has a unique maximum point $(\theta_1^*, \dots, \theta_d^*, \theta_{d+1}^*)$ on the restricted (convex) set (A.13).

The relation between $(\theta, \theta_{d+1})$ and $(\theta, \pi)$ given by (A.11) would give now the only (if any) maximum point of $\mathbb{P} L^s_{\theta,\pi/A} = \mathbb{P} \log g_{\theta,\pi}$ under $\pi \geq 0$.

For the proof of the statements related to the MLE(t), note that by resorting to the canonical form of the exponential family and absorbing $I_A$ into the measure $\lambda$, from Lemma A.2, it is straightforward that the function

$$\log C(\theta) + \sum_{j=1}^d \theta_j \mathbb{P} T_j - \log \mathbb{P}_\theta(A)$$

is a strictly concave function of $\theta$, thus, the results are immediate. □

The next lemma is easily deduced from Theorem 1 in [6].

LEMMA A.3. *Let* $(\mu_0, \Sigma_0) \in \Theta$. *Given* $r_0 > 0$, *let*

$$\mathbb{E}_{(\mu_0,\Sigma_0)}(r_0) := \{\mathcal{E}(\mu, \Sigma, r) : \mathbb{P}_{(\mu_0,\Sigma_0)}[\mathcal{E}(\mu, \Sigma, r)] \geq \mathbb{P}_{(\mu_0,\Sigma_0)}[\mathcal{E}(\mu_0, \Sigma_0, r_0)]\}.$$

*Then, the volume of* $\mathcal{E}(\mu_0, \Sigma_0, r_0)$ *is minimal in* $\mathbb{E}_{(\mu_0,\Sigma_0)}(r_0)$.

PROOF OF PROPOSITION 2.10. Because of the regularity of the model, the map $\theta :\to \mathbb{P}_\theta(\mathcal{E}((\mu_0, \Sigma_0), r))$ is continuously differentiable in a neighborhood of $\theta_0$, so it suffices to show that for every fixed value of $\varsigma^2$ $(= |\Sigma|^{1/p})$ the function has a local maximum at $\mu = \mu_0$ and $\Xi = \Sigma_0/\varsigma_0^2$.

Let $\mu_1$ and $\Sigma_1$ with $|\Sigma_1| = |\Sigma_0|$. Because of the elliptical character of the model, we have that $\mathbb{P}_{(\mu_0,\Sigma_0)}(\mathcal{E}((\mu_0, \Sigma_0), r)) = \mathbb{P}_{(\mu_1,\Sigma_1)}(\mathcal{E}((\mu_1, \Sigma_1), r))$. Thus, if we assume that

$$\mathbb{P}_{(\mu_0,\Sigma_0)}(\mathcal{E}((\mu_0, \Sigma_0), r)) < \mathbb{P}_{(\mu_1,\Sigma_1)}(\mathcal{E}((\mu_0, \Sigma_0), r)),$$

then the absolute continuity of $\mathbb{P}_{(\mu_1,\Sigma_1)}$ implies that, for some $r^* < r$,

$$\mathbb{P}_{(\mu_1,\Sigma_1)}(\mathcal{E}((\mu_0, \Sigma_0), r^*)) = \mathbb{P}_{(\mu_1,\Sigma_1)}(\mathcal{E}((\mu_1, \Sigma_1), r)) = p_0.$$

Then, the volume of the ellipsoid $\mathcal{E}((\mu_0, \Sigma_0), r^*)$ would be strictly lower than $\mathcal{E}((\mu_1, \Sigma_1), r)$ with the same probability, contradicting Lemma A.3. □

Lemmas A.4 and A.5 include some well-known properties, and are stated for reference.

LEMMA A.4. *If* $\mathbb{P}$ *belongs to the GEM given by an elliptical family and* $A$ *is the MVE of* $\mathbb{P}$, *then:*

1. $\mathbb{P}(A) = 1/2$, *and,*
2. *if* $A = \mathcal{E}(\mu, \Sigma, r)$, *then* $\lim_{\varepsilon \to 0+} \mathbb{P}(\mathcal{E}(\mu, \Sigma, r + \varepsilon)) = 1/2$.

The next lemma follows from the well-known fact that the class

$$\mathcal{C} := \{\{x \in \mathbb{R}^p : |\langle x - \mu, v\rangle| \leq d\} : \mu \in \mathbb{R}^p, v \in S(0,1) \text{ and } d > 0\},$$

and the class of all ellipsoids constitute two Vapnik–Cervonenkis (VC) classes.

LEMMA A.5. *Let* $\{X_n\}_n$ *be a random sample taken from a probability distribution* $\mathbb{P}$, *then:*

1. $\sup\{|P_n(A) - \mathbb{P}(A)| : A \textit{ is an ellipsoid }\} \to 0$, *a.s.*
2. $\sup\{|P_n(A \cap B) - \mathbb{P}(A \cap B)| : A \textit{ is an ellipsoid and } B \in \mathcal{C}\} \to 0$, *a.s.*

LEMMA A.6. *Let $g$ be a decreasing function which defines an elliptical family on $\mathbb{R}^p$. Let $\mathbb{P} = (1 - \pi_0)\mathbb{P}_{\theta_0} + \pi_0 Q$, $0 \le \pi_0 < 1/2$ and $\theta_0 \in \Theta$, be a distribution in the GEM of the elliptical family defined by $g$.*

*Let $A = \mathcal{E}(\mu_0, \Sigma_0, r_0)$ be the MVE of $\mathbb{P}$. Then, for every $\eta > 0$, there exist $d > 0$ and $\varepsilon > 0$, such that, if we denote $A^\varepsilon = \mathcal{E}(\mu, \Sigma, r + \varepsilon)$, then*

$$\sup_{v \in S(0,1)} \sup_{\mu \in \mathbb{R}^p} \mathbb{P}[A^\varepsilon \cap \{x \in \mathbb{R}^p : |\langle x - \mu, v\rangle| \le d\}] < \eta.$$

PROOF. Let $\eta > 0$. Since $\mathbb{P}_{\theta_0}$ is absolutely continuous, it is easy to show that there exists $d > 0$ such that

$$\sup_{v \in S(0,1)} \sup_{\mu \in \mathbb{R}^p} \mathbb{P}[\{x \in \mathbb{R}^p : |\langle x - \mu, v\rangle| \le d\}] < \eta.$$

This and (2) in Lemma A.4 give the result. $\square$

LEMMA A.7. *Assume the hypotheses of Theorem 3.2. Let $\delta_n^1, \dots, \delta_n^p$ be the eigenvalues of $\hat{\Sigma}_n$. Let $\Delta_n = \inf\{\delta_n^1, \dots, \delta_n^p\}$. Then*

$$\liminf_n \Delta_n > 0, \qquad a.s.$$

PROOF. We will only treat the case $p > 1$. Let $\gamma > p/2$ be the value associated by Gp to $g$. Obviously, $2\gamma - p > 0$ and, then, there exists $\eta > 0$ such that $2\gamma(1 - 2\eta) > p$.

From Lemma A.6, there exist $\varepsilon > 0$ and $d > 0$ such that

$$\sup_{v \in S(0,1)} \sup_{\mu \in \mathbb{R}^p} \mathbb{P}[A^\varepsilon \cap \{x \in \mathbb{R}^p : |\langle x - \mu, v\rangle| \le d\}] < \eta.$$

Taking into account Proposition 3.1, Lemmas A.5 and A.4 and that $I_{A^\varepsilon} \log(f_{\theta_0})$ is $\mathbb{P}$-integrable, we have that there exists a probability one set $\Omega_0$ such that if $\omega \in \Omega_0$, then

$$P_n^\omega[A_n] \to \mathbb{P}(A) = 1/2. \tag{A.14}$$

There exists $N \in \mathbb{N}$ ($N = N(\omega)$) such that if $n \ge N$, then $A_n \subset A^\varepsilon$ and,

$$\sup_{v \in S(0,1)} \sup_{\mu \in \mathbb{R}^p} \mathbb{P}[A_n \cap \{x \in \mathbb{R}^p : |\langle x - \mu, v\rangle| \le d\}] < \eta, \tag{A.15}$$

and

$$P_n^\omega I_{A_n} \log(f_{\theta_0}) \to \mathbb{P} I_A \log(f_{\theta_0}). \tag{A.16}$$

Let $\omega \in \Omega_0$. $\omega$ will remain fixed and will be omitted in the notation.

Statement (A.14) implies that requirements on $m$ in Proposition 2.5 hold from an index onward, and then the MLE(s) exists from this index onward.

Let us denote $\hat{\mu}_n = (\mu_n^1, \dots, \mu_n^p)$. Let $j_n$ be such that $\Delta_n = \delta_n^{j_n}$. Let us assume that there exists a subsequence such that $\lim_k \Delta_{n_k} = 0$. To simplify the notation, we will denote this subsequence with the same notation as the original one. Let us consider the set

$$B_n := \{x = (x^1, \ldots, x^p) \in \mathbb{R}^p : (x^{j_n} - \mu_n^{j_n})^2 \geq d^2\}.$$

From (A.15) we have that $P_n[B_n \cap A_n] \geq P_n(A_n) - \eta$ and now the proof goes by repeating the same steps as in Proposition 2.4 with the set $I_k$ being replaced there by the set $B_n$ here, because we have

$$\begin{aligned}
&\prod_{i=1}^{n} (f_{(\hat{\mu}_n, \hat{\Sigma}_n)}(x_i))^{I_{A_n}(x_i)} \\
&\quad \leq (\Delta_n)^{-(p/2)nP_n(A_n)} g(0)^{n\eta} d^{-2\gamma n(P_n(A_n)-\eta)} (\Delta_n)^{\gamma n(P_n(A_n)-\eta)} \qquad \text{(A.17)}\\
&\quad = ((\Delta_n)^{2\gamma(P_n(A_n)-\eta)-pP_n(A_n)/2} g(0)^{\eta} d^{-2\gamma(P_n(A_n)-\eta)})^n. \qquad \square
\end{aligned}$$

Proof of Theorem 3.2. Let us denote $\hat{\theta}_{s,n} = (\hat{\pi}_n, \hat{\mu}_n, \hat{\Sigma}_n)$ and let $A_n$ be an empirical MVE. First, we will prove that the sequence $\{(\hat{\pi}_n, \hat{\mu}_n, \hat{\Sigma}_n)\}_n$ is a.s. included in a compact subset of $[0,1) \times \Theta$. To this, let $\Omega_0$ be a probability one set whose points satisfy 1 in Lemma A.5, Proposition 3.1, (A.14), (A.16) and Lemma A.7, and let $\omega \in \Omega_0$. This point will remain fixed throughout the proof and we will make no reference to it in the notation.

The second term in (2.8) does not depend on $\pi$ and the third one is bounded. Since $(\hat{\pi}_n, \hat{\mu}_n, \hat{\Sigma}_n)$ maximizes (2.8) and the first term converges to $-\infty$ if $\hat{\pi}_n \to 1$, it may not happen that $\limsup \hat{\pi}_n = 1$.

Let $\delta_n^1, \ldots, \delta_n^p$ be the eigenvalues of $\hat{\Sigma}_n$. Let $\Delta^n = \sup(\delta_n^1, \ldots, \delta_n^p)$. Following the same steps as in Proposition 2.4, we would prove that if it were $\limsup \Delta^n = \infty$, then there would exist a subsequence such that

$$\lim_k \prod_{i=1}^{n_k} (f_{(\hat{\mu}_{n_k}, \hat{\Sigma}_{n_k})}(x_i))^{I_{A_{n_k}}(x_i)} = 0$$

and, for this sequence, the second term in (2.8) would converge to $-\infty$, which is impossible because (A.16) is satisfied.

Therefore, the sequence $\{(\hat{\pi}_n, \hat{\mu}_n, \hat{\Sigma}_n)\}_n$ is a.s. included in a compact subset of $[0,1) \times \Theta$. Let us consider a convergent subsequence $\{(\hat{\pi}_{n_k}, \hat{\mu}_{n_k}, \hat{\Sigma}_{n_k})\}_k$ with limit $(\pi^*, \mu^*, \Sigma^*)$.

Lemma A.4 and (A.14) trivially give that

$$P_{n_k}(A_{n_k}) \log(1 - \pi_{n_k}) \to \mathbb{P}(A) \log(1 - \pi^*). \qquad \text{(A.18)}$$

Proposition 3.1 implies that if we denote $A = \mathcal{E}(\mu, \Sigma, r)$, then for every $\varepsilon > 0$ there exists $N$ such that, if $k \geq N$, then

$$\mathcal{E}(\mu, \Sigma, r - \varepsilon) \subset A_{n_k} \subset \mathcal{E}(\mu, \Sigma, r + \varepsilon). \qquad \text{(A.19)}$$

However, if we denote $\hat{\theta}_{n_k} = (\hat{\mu}_{n_k}, \hat{\Sigma}_{n_k})$, $\theta^* = (\mu^*, \Sigma^*)$, (2.12) implies that

$$\lim_k \mathbb{P}_{\hat{\theta}_{n_k}}(\mathcal{E}(\mu, \Sigma, r^\varepsilon)) = \mathbb{P}_{\theta^*}(\mathcal{E}(\mu, \Sigma, r^\varepsilon)),$$

where $r^\varepsilon \in \{r - \varepsilon, r + \varepsilon\}$. From here, (A.19) and the continuity of $\mathbb{P}_{\theta^*}$, we obtain that

$$\lim_k \mathbb{P}_{\hat{\theta}_{n_k}}(A^c_{n_k}) = \mathbb{P}_{\theta^*}(A^c).$$

This, with (A.14), gives that

$$\begin{aligned} &\lim_k P_{n_k}(A^c_{n_k}) \log((1 - \pi_{n_k})\mathbb{P}_{\hat{\theta}_{n_k}}(A^c_{n_k}) + \pi_{n_k}) \\ &\quad = \mathbb{P}(A^c) \log((1 - \pi^*)\mathbb{P}_{\theta^*}(A^c) + \pi^*). \end{aligned} \tag{A.20}$$

The continuity of $g$ implies that $\lim_k f_{\hat{\theta}_{n_k}} = f_{\theta^*}$. Moreover, the sequence $\{f_{\hat{\theta}_{n_k}}\}_k$ is uniformly bounded by a constant because the sequence $\{\hat{\Sigma}_{n_k}\}_k$ is contained in a compact subset of $\Theta$ and $g$ is bounded by $g(0)$. Thus, taking into account that 1 in Lemma A.5 implies that the sequence of distributions $\{P_{n_k}\}_k$ converges in distribution to $\mathbb{P}$, and that $A$ is a continuity set of $\mathbb{P}$, it is a standard exercise to prove that

$$P_{n_k}(I_{A_{n_k}} \log(f_{\theta_{n_k}})) \to \mathbb{P}(I_A \log(f_{\theta^*})).$$

This, (A.18) and (A.20) give that

$$\lim_k P_{n_k} L^s_{\theta_{n_k}, \pi_{n_k}/A_{n_k}} = \mathbb{P} L^s_{\theta^*, \pi^*/A}. \tag{A.21}$$

On the other hand, from the assumptions on $\Omega_0$, it can be deduced that

$$\lim_k P_{n_k} L^s_{\theta_0, \pi_0/A_{n_k}} = \mathbb{P} L^s_{\theta_0, \pi_0/A}. \tag{A.22}$$

But, from Proposition 2.3, we obtain that $\mathbb{P} L^s_{\theta^*, \pi^*/A} \le \mathbb{P} L^s_{\theta_0, \pi_0/A}$, and, by definition of the smart estimate, we also have that

$$\lim_k P_{n_k} L^s_{\theta_0, \pi_0/A_{n_k}} \le \lim_k P_{n_k} L^s_{\theta_{n_k}, \pi_{n_k}/A_{n_k}}.$$

This, (A.21), (A.22) and the inequality (2.11) imply that $\theta_0 = \theta^*$.

The consistency of the MLE(c) under the elliptical model can be proved with the same scheme by considering the easier case $\pi = 0$.

The proof for the MLE(r) follows the same steps, once we show that the dependence on $\alpha$ of the restrictions given by $\Theta^\alpha_A$ does not constitute any constraint from the asymptotic point of view. This is proved in Proposition 3.2. □

Proof of Proposition 3.2. Let $\delta \in (0, 1/2 - \alpha)$. From Proposition 3.1 and the Glivenko–Cantelli property of the class of ellipsoids, we deduce that, for $n \geq n_0$ large enough (and depending on the sample), $\mathbb{P}_{\theta_0}(A_n) > \alpha + \eta$ holds. On the other hand, (2.12) shows that for $\varepsilon > 0$ there exists $\delta > 0$ such that $\sup_{B \in \beta^p} |\mathbb{P}_\theta(B) - \mathbb{P}_{\theta_0}(B)| < \varepsilon$ whenever $\|\theta - \theta_0\| < \delta$. Both relations give that $\mathbb{P}_\theta(A_n) > \alpha$, so that $\theta \in \Theta^\alpha_{A_n}$, if $\|\theta - \theta_0\| < \delta$ and $n > n_0$. □

Now we will adapt Section 3.2.4 in [26] to our semiparametric setup. We only include some keys for the adapted proofs which verbatim would repeat the arguments there.

THEOREM A.1 (Extension of Theorem 3.2.16 in [26]). *Let $\{M_n\}_n$ be stochastic processes, all of them indexed by the same product $\Theta \times K$ of an open subset $\Theta$ and a compact subset $K$ of two Euclidean spaces, and $M : \Theta \times K \to \mathbb{R}$ be a deterministic function.*

*Assume that for every $\gamma \in K$ the function $\theta \to M(\theta, \gamma)$ has a unique maximum $\theta_0$ where it is twice continuously differentiable w.r.t. $\theta$, with nonsingular continuous (w.r.t. $\gamma$) second derivative matrix $V(\gamma)$. Suppose also that*

$$\begin{aligned}\sqrt{n}(M_n(\theta_n, \gamma_n) - M(\theta_n, \gamma_n)) &- \sqrt{n}(M_n(\theta_0, \gamma_n) - M(\theta_0, \gamma_n)) \\ &= (\theta_n - \theta_0)^T Z_n(\theta_n, \gamma_n) + o^*_P(\|\theta_n - \theta_0\|^2) \\ &\quad + o^*_P(\|\theta_n - \theta_0\| + \sqrt{n}\|\theta_n - \theta_0\|^2 + n^{-1/2})\end{aligned}$$

*for every sequence $\theta_n = \theta_0 + o^*_P(1)$, every sequence $\{\gamma_n\} \subset K$, and a uniformly tight sequence $Z_n(\theta_n, \gamma_n)$ of random vectors.*

*If the sequence $\hat{\theta}_n(\gamma_n)$ converges in outer probability to $\theta_0$ and satisfies*

$$M_n(\hat{\theta}_n(\gamma_n), \gamma_n) \geq \sup_\theta M_n(\theta, \gamma_n) - o_P(n^{-1}),$$

*for every $n \in \mathbb{N}$, then $\sqrt{n}(\hat{\theta}_n(\gamma_n) - \theta_0) = -(V(\gamma_n))^{-1} Z_n(\hat{\theta}_n(\gamma_n), \gamma_n) + o^*_P(1)$.*

PROOF. For every sequence $h_n = o^*_P(1)$, the hypotheses yield

$$\begin{aligned}M_n(\theta_0 + h_n, \gamma_n) &- M_n(\theta_0, \gamma_n) \\ &= \tfrac{1}{2} h'_n V(\gamma_n) h_n + n^{-1/2} h'_n Z_n(\theta_0 + h_n, \gamma_n) \\ &\quad + o^*_P(\|h_n\|^2 + (\sqrt{n})^{-1}\|h_n\| + n^{-1}).\end{aligned} \tag{A.23}$$

Take $h_n = \widehat{\theta}_n(\gamma_n) - \theta_0$, and follow the proof in [26], taking into account that the term $h'_n V(\gamma_n) h_n$ on the right-hand side is bounded above by $c\|h_n\|^2$ for some $c > 0$. This holds because, on the contrary, there should exist a sequence $\{\gamma_n\} \subset K$ such that the corresponding sequence of minimum eigenvalues of $V(\gamma_n)$ would converge to 0. Then, for some convergent subsequence to some $\gamma \in K$, the continuity of $V$ would give the contradiction of the singularity of $V(\gamma)$. [The same argument based on the compactness of $K$ makes it possible to guarantee that the eigenvalues of $V(\gamma_n)$ are bounded, so that $-n^{-1/2}(V(\gamma_n))^{-1}Z_n(\hat{\theta}_n(\gamma_n),\gamma_n)$ is $O_P^*(n^{-1/2})$, and, hence, to apply (A.23) also to $h_n = -n^{-1/2}(V(\gamma_n))^{-1}Z_n(\hat{\theta}_n(\gamma_n),\gamma_n)$ analogously to the original proof.] $\square$

Now let $\Theta$ be the parameter space in the elliptical model and let $K$ be a compact subset of $\Gamma$. Let $m_{\theta,\gamma}$ be a real function, consider the sample probability, $P_n$, corresponding to $n$ i.i.d. observations from $\mathbb{P}$, $M_n(\theta,\gamma) = P_n m_{\theta,\gamma}$ and $M(\theta,\gamma) = \mathbb{P}m_{\theta,\gamma}$, as well as the empirical process $\mathbb{G}_n m_{\theta\gamma} = \sqrt{n}(M_n(\theta,\gamma) - M(\theta,\gamma))$. The differentiability involved in the preceding theorem can be guaranteed through the condition required in Lemma A.8.

LEMMA A.8 (Extension of Lemma 3.2.19 in [26]). *Suppose that for every $\gamma$ in the compact set $K$, there exists a vector-valued function $\dot{m}_{\theta_0,\gamma}$ such that, for some $\delta > 0$,*

$$\left\{\frac{m_{\theta\gamma} - m_{\theta\gamma_0} - (\theta-\theta_0)^T \dot{m}_{\gamma\theta_0}}{\|\theta-\theta_0\|} : \|\theta-\theta_0\| < \delta,\ \gamma \in K\right\}$$

*is $\mathbb{P}$-Donsker and that, uniformly for $\gamma \in K$,*

$$\mathbb{P}(m_{\theta\gamma} - m_{\theta_0\gamma} - (\theta-\theta_0)^T \dot{m}_{\theta\gamma_0})^2 = o(\|\theta-\theta_0\|^2).$$

*Then, if $\theta_n = \theta_0 + o_P^*(1)$, we have that, uniformly in $\gamma$,*

$$\begin{aligned}\mathbb{G}_n(m_{\theta_n\gamma} - m_{\theta_0\gamma}) \\ = (\theta_n-\theta_0)^T \mathbb{G}_n \dot{m}_{\theta_0\gamma} + o_P^*(\|\theta_n-\theta_0\| + \sqrt{n}\|\theta_n-\theta_0\|^2 + n^{-1/2}).\end{aligned} \tag{A.24}$$

PROOF. It suffices to adapt the proof in [26] to the function

$$f : \ell^\infty(\Theta_d \times K) \times (\Theta_d \times K) \to \mathbb{R}^d$$

given by $f(z,(\theta,\gamma)) = z(\theta,\gamma)$ ($\Theta_d := \{\|\theta-\theta_0\| < \delta\}$) and the stochastic processes

$$Z_n(\theta,\gamma) = \mathbb{G}_n \frac{m_{\theta\gamma} - m_{\theta_0\gamma} - (\theta-\theta_0)^T \dot{m}_{\theta_0\gamma}}{\|\theta-\theta_0\|},$$

which, from the hypotheses, converge in $\ell^\infty(\Theta_d \times K)$ to a tight Gaussian process $Z$. $\square$

PROOF OF LEMMA 3.1. In order to simplify the notation, given $x \in \mathbb{R}^p$ and $\theta = (\mu,\Sigma)$, let us denote $x_\theta = (x-\mu)^T\Sigma^{-1}(x-\mu)$.

From the continuity and the nonincreasing character of $g$, we deduce that there exist $K$, a compact neighborhood of $\gamma_0$, $\delta > 0$, and an ellipsoid $\mathcal{E}(\gamma^*)$ such that $\bigcup_{\gamma\in K}\mathcal{E}(\gamma)\subset\mathcal{E}(\gamma^*)$ and, if $V_\delta:=\{\theta\in\Theta:\|\theta-\theta_0\|<\delta\}$, then $\inf_{x\in\mathcal{E}(\gamma^*),\theta\in V_\delta} g(x_{\theta_0})>0$.

Since the set $\{x_\theta:\theta\in\Theta_\delta, x\in\mathcal{E}(\gamma^*)\}$ is bounded, the second statement follows from the derivability w.r.t. $\theta$, leading to the Fréchet derivability in $L^2$. Also, note that the hypothesis on the continuity of the second derivative implies that $g'$ is Lipschitz in its effective domain.

Thus, the components of $\dot{m}_{\theta\gamma}$ are easily seen to be

$$\frac{\int_{I_{\mathcal{E}(\gamma)}} 2\Sigma^{-1}(y-\mu)g'(y_\theta)\,dy}{\int_{I_{\mathcal{E}(\gamma)}} g(y_\theta)\,dy}-\frac{2\Sigma^{-1}(x-\mu)g'(x_\theta)}{g(x_\theta)}$$

for the derivative w.r.t. $\mu$, while those corresponding to $\Sigma$ are

$$\frac{\int_{I_{\mathcal{E}(\gamma)}} 2\Sigma^{-1}(y-\mu)(y-\mu)^T\Sigma^{-1}g'(y_\theta)\,dy}{\int_{I_{\mathcal{E}(\gamma)}} g(y_\theta)\,dy}-\frac{2\Sigma^{-1}(x-\mu)(x-\mu)^T\Sigma^{-1}g'(x_\theta)}{g(x_\theta)}.$$

To verify that (3.5) is $\mathbb{P}$-Donsker, we only need to give some steps concerning the permanence of the Donsker property as developed in Section 2.10 in

TABLE 1

*Asymptotic efficiencies to estimate an element of $\mu$*

| | | **Dimension** | | | | |
|---|---|---|---|---|---|---|
| | | $p=2$ | $p=3$ | $p=5$ | $p=10$ | $p=30$ |
| Gaussian | MLE(c) | 0.1531 | 0.2032 | 0.2613 | 0.3263 | 0.3984 |
| | $\text{MLE(c)}_{0.25}$ | 0.4049 | 0.4658 | 0.5301 | 0.5991 | 0.6627 |
| | $\text{MLE(c)}_{0.10}$ | 0.6675 | 0.7184 | 0.7650 | 0.8085 | 0.8503 |
| | $\text{MLE(c)}_{0.025}$ | 0.8821 | 0.9040 | 0.9242 | 0.9414 | 0.9579 |
| $t_1$ | MLE(c) | 0.5147 | 0.5933 | 0.6414 | 0.6434 | 0.5975 |
| | $\text{MLE(c)}_{0.25}$ | 0.8037 | 0.8374 | 0.8542 | 0.8492 | 0.8198 |
| | $\text{MLE(c)}_{0.10}$ | 0.9387 | 0.9482 | 0.9530 | 0.9490 | 0.9350 |
| | $\text{MLE(c)}_{0.025}$ | 0.9884 | 0.9900 | 0.9906 | 0.9896 | 0.9855 |
| $t_5$ | MLE(c) | 0.2889 | 0.3780 | 0.4713 | 0.5481 | 0.5737 |
| | $\text{MLE(c)}_{0.25}$ | 0.6132 | 0.6839 | 0.7465 | 0.7931 | 0.8047 |
| | $\text{MLE(c)}_{0.10}$ | 0.8450 | 0.8792 | 0.9070 | 0.9249 | 0.9280 |
| | $\text{MLE(c)}_{0.025}$ | 0.9636 | 0.9724 | 0.9793 | 0.9833 | 0.9839 |
| $t_8$ | MLE(c) | 0.2469 | 0.3296 | 0.4193 | 0.5092 | 0.5583 |
| | $\text{MLE(c)}_{0.25}$ | 0.5586 | 0.6339 | 0.7067 | 0.7657 | 0.7953 |
| | $\text{MLE(c)}_{0.10}$ | 0.8083 | 0.8495 | 0.8850 | 0.9133 | 0.9240 |
| | $\text{MLE(c)}_{0.025}$ | 0.9511 | 0.9631 | 0.9729 | 0.9800 | 0.9826 |
| $t_{15}$ | MLE(c) | 0.2082 | 0.2805 | 0.3640 | 0.4544 | 0.5306 |
| | $\text{MLE(c)}_{0.25}$ | 0.5025 | 0.5769 | 0.6535 | 0.7249 | 0.7768 |
| | $\text{MLE(c)}_{0.10}$ | 0.7642 | 0.8108 | 0.8535 | 0.8912 | 0.9157 |
| | $\text{MLE(c)}_{0.025}$ | 0.9337 | 0.9486 | 0.9627 | 0.9734 | 0.9801 |

Table 2
*Asymptotic efficiencies to estimate a diagonal element of* $\Sigma$

| | | Dimension | | | | |
|---|---|---|---|---|---|---|
| | | $p = 2$ | $p = 3$ | $p = 5$ | $p = 10$ | $p = 30$ |
| Gaussian | MLE(t) | 0.0266 | 0.0521 | 0.1023 | 0.1790 | 0.3000 |
| | MLE$(t)_{0.25}$ | 0.1375 | 0.2059 | 0.2955 | 0.4184 | 0.5593 |
| | MLE$(t)_{0.10}$ | 0.3594 | 0.4457 | 0.5508 | 0.6657 | 0.7744 |
| | MLE$(t)_{0.025}$ | 0.6673 | 0.7364 | 0.8072 | 0.8674 | 0.9273 |
| | MLE(c) | 0.2666 | 0.2293 | 0.2161 | 0.2392 | 0.3206 |
| | MLE$(c)_{0.25}$ | 0.4560 | 0.4217 | 0.4248 | 0.4813 | 0.5793 |
| | MLE$(c)_{0.10}$ | 0.6551 | 0.6321 | 0.6534 | 0.7128 | 0.7894 |
| | MLE$(c)_{0.025}$ | 0.8408 | 0.8435 | 0.8614 | 0.8918 | 0.9336 |
| $t_1$ | MLE(t) | 0.2004 | 0.2990 | 0.3941 | 0.4597 | 0.4938 |
| | MLE$(t)_{0.25}$ | 0.4941 | 0.5976 | 0.6778 | 0.7244 | 0.7457 |
| | MLE$(t)_{0.10}$ | 0.7351 | 0.8126 | 0.8599 | 0.8879 | 0.8968 |
| | MLE$(t)_{0.025}$ | 0.8778 | 0.9334 | 0.9619 | 0.9712 | 0.9747 |
| | MLE(c) | 0.4255 | 0.3736 | 0.4085 | 0.4611 | 0.4938 |
| | MLE$(c)_{0.25}$ | 0.6619 | 0.6507 | 0.6873 | 0.7251 | 0.7458 |
| | MLE$(c)_{0.10}$ | 0.8486 | 0.8480 | 0.8667 | 0.8884 | 0.8968 |
| | MLE$(c)_{0.025}$ | 0.9588 | 0.9593 | 0.9661 | 0.9715 | 0.9747 |
| $t_5$ | MLE(t) | 0.0786 | 0.1492 | 0.2512 | 0.3749 | 0.4664 |
| | MLE$(t)_{0.25}$ | 0.3028 | 0.4134 | 0.5381 | 0.6518 | 0.7279 |
| | MLE$(t)_{0.10}$ | 0.5914 | 0.6877 | 0.7773 | 0.8498 | 0.8903 |
| | MLE$(t)_{0.025}$ | 0.8415 | 0.8890 | 0.9282 | 0.9571 | 0.9690 |
| | MLE(c) | 0.3661 | 0.3064 | 0.3118 | 0.3865 | 0.4670 |
| | MLE$(c)_{0.25}$ | 0.5749 | 0.5474 | 0.5854 | 0.6607 | 0.7283 |
| | MLE$(c)_{0.10}$ | 0.7736 | 0.7739 | 0.8074 | 0.8549 | 0.8906 |
| | MLE$(c)_{0.025}$ | 0.9269 | 0.9301 | 0.9433 | 0.9597 | 0.9692 |
| $t_8$ | MLE(t) | 0.0609 | 0.1182 | 0.2116 | 0.3347 | 0.4502 |
| | MLE$(t)_{0.25}$ | 0.2552 | 0.3611 | 0.4883 | 0.6199 | 0.7199 |
| | MLE$(t)_{0.10}$ | 0.5411 | 0.6430 | 0.7429 | 0.8248 | 0.8847 |
| | MLE$(t)_{0.025}$ | 0.8195 | 0.8701 | 0.9144 | 0.9500 | 0.9683 |
| | MLE(c) | 0.3437 | 0.2881 | 0.2874 | 0.3545 | 0.4514 |
| | MLE$(c)_{0.25}$ | 0.5481 | 0.5175 | 0.5517 | 0.6347 | 0.7207 |
| | MLE$(c)_{0.10}$ | 0.7488 | 0.7465 | 0.7816 | 0.8342 | 0.8852 |
| | MLE$(c)_{0.025}$ | 0.9131 | 0.9173 | 0.9327 | 0.9541 | 0.9685 |
| $t_{15}$ | MLE(t) | 0.0458 | 0.0918 | 0.1717 | 0.2908 | 0.4302 |
| | MLE$(t)_{0.25}$ | 0.2089 | 0.3038 | 0.4288 | 0.5684 | 0.6981 |
| | MLE$(t)_{0.10}$ | 0.4796 | 0.5785 | 0.6890 | 0.7962 | 0.8735 |
| | MLE$(t)_{0.025}$ | 0.7771 | 0.8400 | 0.8961 | 0.9364 | 0.9669 |
| | MLE(c) | 0.3168 | 0.2685 | 0.2630 | 0.3205 | 0.4331 |
| | MLE$(c)_{0.25}$ | 0.5177 | 0.4855 | 0.5118 | 0.5931 | 0.7003 |
| | MLE$(c)_{0.10}$ | 0.7188 | 0.7098 | 0.7452 | 0.8106 | 0.8746 |
| | MLE$(c)_{0.025}$ | 0.8931 | 0.8988 | 0.9189 | 0.9433 | 0.9673 |

Table 3
*Asymptotic efficiencies to estimate an off-diagonal element of* $\Sigma$

| | | Dimension | | | | |
|---|---|---|---|---|---|---|
| | | $p = 2$ | $p = 3$ | $p = 5$ | $p = 10$ | $p = 30$ |
| Gaussian | MLE(t) | 0.0332 | 0.0631 | 0.1130 | 0.1929 | 0.3030 |
| | $\text{MLE(t)}_{0.25}$ | 0.1621 | 0.2323 | 0.3247 | 0.4361 | 0.5718 |
| | $\text{MLE(t)}_{0.10}$ | 0.4070 | 0.4874 | 0.5854 | 0.6890 | 0.7872 |
| | $\text{MLE(t)}_{0.025}$ | 0.7151 | 0.7716 | 0.8333 | 0.8805 | 0.9315 |
| | MLE(c) | 0.0332 | 0.0631 | 0.1130 | 0.1929 | 0.3030 |
| | $\text{MLE(c)}_{0.25}$ | 0.1621 | 0.2323 | 0.3247 | 0.4361 | 0.5718 |
| | $\text{MLE(c)}_{0.10}$ | 0.4070 | 0.4874 | 0.5854 | 0.6890 | 0.7872 |
| | $\text{MLE(c)}_{0.025}$ | 0.7151 | 0.7716 | 0.8333 | 0.8805 | 0.9315 |
| $t_1$ | MLE(t) | 0.0581 | 0.0997 | 0.1540 | 0.2064 | 0.2371 |
| | $\text{MLE(t)}_{0.25}$ | 0.1470 | 0.2023 | 0.2655 | 0.3219 | 0.3598 |
| | $\text{MLE(t)}_{0.10}$ | 0.2202 | 0.2757 | 0.3390 | 0.3962 | 0.4331 |
| | $\text{MLE(t)}_{0.025}$ | 0.2646 | 0.3178 | 0.3784 | 0.4319 | 0.4701 |
| | MLE(c) | 0.7773 | 0.7684 | 0.7634 | 0.7598 | 0.7553 |
| | $\text{MLE(c)}_{0.25}$ | 0.8665 | 0.8697 | 0.8726 | 0.8762 | 0.8782 |
| | $\text{MLE(c)}_{0.10}$ | 0.9401 | 0.9434 | 0.9468 | 0.9497 | 0.9505 |
| | $\text{MLE(c)}_{0.025}$ | 0.9836 | 0.9852 | 0.9861 | 0.9871 | 0.9876 |
| $t_5$ | MLE(t) | 0.0301 | 0.0572 | 0.1021 | 0.1643 | 0.2244 |
| | $\text{MLE(t)}_{0.25}$ | 0.1119 | 0.1584 | 0.2173 | 0.2878 | 0.3479 |
| | $\text{MLE(t)}_{0.10}$ | 0.2164 | 0.2598 | 0.3122 | 0.3743 | 0.4271 |
| | $\text{MLE(t)}_{0.025}$ | 0.3003 | 0.3332 | 0.3721 | 0.4217 | 0.4659 |
| | MLE(c) | 0.6862 | 0.6922 | 0.7047 | 0.7263 | 0.7444 |
| | $\text{MLE(c)}_{0.25}$ | 0.7691 | 0.7927 | 0.8198 | 0.8495 | 0.8687 |
| | $\text{MLE(c)}_{0.10}$ | 0.8728 | 0.8950 | 0.9166 | 0.9341 | 0.9472 |
| | $\text{MLE(c)}_{0.025}$ | 0.9578 | 0.9674 | 0.9756 | 0.9820 | 0.9866 |
| $t_8$ | MLE(t) | 0.0261 | 0.0493 | 0.0899 | 0.1497 | 0.2194 |
| | $\text{MLE(t)}_{0.25}$ | 0.1047 | 0.1488 | 0.2053 | 0.2753 | 0.3420 |
| | $\text{MLE(t)}_{0.10}$ | 0.2166 | 0.2594 | 0.3074 | 0.3659 | 0.4215 |
| | $\text{MLE(t)}_{0.025}$ | 0.3181 | 0.3451 | 0.3788 | 0.4206 | 0.4648 |
| | MLE(c) | 0.6518 | 0.6625 | 0.6809 | 0.7093 | 0.7387 |
| | $\text{MLE(c)}_{0.25}$ | 0.7311 | 0.7610 | 0.7966 | 0.8342 | 0.8631 |
| | $\text{MLE(c)}_{0.10}$ | 0.8422 | 0.8718 | 0.9000 | 0.9250 | 0.9426 |
| | $\text{MLE(c)}_{0.025}$ | 0.9446 | 0.9571 | 0.9695 | 0.9791 | 0.9854 |
| $t_{15}$ | MLE(t) | 0.0220 | 0.0424 | 0.0775 | 0.1317 | 0.2068 |
| | $\text{MLE(t)}_{0.25}$ | 0.0967 | 0.1379 | 0.1913 | 0.2589 | 0.3321 |
| | $\text{MLE(t)}_{0.10}$ | 0.2146 | 0.2564 | 0.3061 | 0.3587 | 0.4157 |
| | $\text{MLE(t)}_{0.025}$ | 0.3348 | 0.3607 | 0.3884 | 0.4230 | 0.4599 |
| | MLE(c) | 0.6097 | 0.6236 | 0.6484 | 0.6827 | 0.7287 |
| | $\text{MLE(c)}_{0.25}$ | 0.6841 | 0.7191 | 0.7617 | 0.8113 | 0.8541 |
| | $\text{MLE(c)}_{0.10}$ | 0.8015 | 0.8371 | 0.8756 | 0.9108 | 0.9383 |
| | $\text{MLE(c)}_{0.025}$ | 0.9227 | 0.9413 | 0.9589 | 0.9728 | 0.9832 |

[26], starting from the Lipschitz property of $g'$, that with Theorem 2.10.6 there leads to the Donsker property of the class $\{g'(x_\theta) : \|\theta - \theta_0\| < \delta\}$.

The uniform (below or upper) bounds on the compact set $\mathcal{E}(\gamma^*)$ containing the ellipsoids in the class permit us to apply the properties in Examples 2.10.8 and 2.10.9 in [26] and conclude Donsker's property of the class in (3.5). □

## APPENDIX B: ASYMPTOTIC EFFICIENCY

Tables 1–3 show the efficiency of the proposed estimators in the estimation of an element of $\mu$, and arbitrary diagonal and off-diagonal elements of $\Sigma$ in several dimensions, for the multivariate Gaussian and some $t$ distributions.

We analyze the MLE(c) and MLE(t) and the estimators based on enlarged versions of the MVE to cover $1 - \alpha$ of the theoretical probability [MLE(c)$_\alpha$ or MLE(t)$_\alpha$]. This assures maximum BP of our equivariant estimators.

When estimating the components of $\mu$, the efficiencies of the truncated and censored estimates coincide, and we only show those of the censured one.

The efficiencies have been computed comparing the values of the asymptotic variances (in Theorem 3.3) with the Cramer–Rao bound. The involved integrals have been computed by the Monte Carlo method with 500,000 repetitions.

**Acknowledgments.** The authors want to thank, sincerely, the comments received from two anonymous referees which have improved considerably several points in the paper.

J. A. Cuesta-Albertos
Departamento Matemáticas,
Estadística y Computación
Universidad de Cantabria
Avda. los Castros s.n.
39005 Santander
Spain
E-mail: cuestaj@unican.es

C. Matrán
A. Mayo-Iscar
Departamento de Estadística
e Investigación Operativa
Universidad de Valladolid
Prado de la Magdalena s.n.
47005 Valladolid
Spain
E-mail: matran@eio.uva.es
agustin@med.uva.es